\documentclass{article}
\usepackage[usenames,dvipsnames]{color}
\usepackage{booktabs}
\usepackage{graphicx}
\usepackage{caption}
\usepackage{amssymb}
\usepackage{amsmath,amsthm}
\usepackage[T1]{fontenc}
\usepackage{eufrak}
\usepackage[cp1250]{inputenc}
\usepackage[polish,english]{babel}
\usepackage{subfigure}

\newcommand{\bR}{{\mathbb{R}}}
\newcommand{\bN}{{\mathbb{N}}}

\newtheorem{Theorem}{Theorem}[section]
\newtheorem{Proposition}[Theorem]{Proposition}
\newtheorem{Definition}[Theorem]{Definition}
\newtheorem{Lemma}[Theorem]{Lemma}
\newtheorem{Remark}[Theorem]{Remark}
\newtheorem{Assumption}[Theorem]{Assumption}

\newenvironment{proofof}
{\smallskip\noindent{\textbf{Proof.-}}
\hspace{1pt}}{\hspace{-5pt}{\nobreak\nobreak\hfill\nobreak
$\square$\vspace{2pt}\par}\smallskip\goodbreak}

\begin{document}
\title{The Escalator Boxcar Train Method for a System of Aged-structured Equations in the Space of Measures}

\author{
Jos\'e A. Carrillo
\footnote{Department of Mathematics, Imperial College London, SW7 2AZ London, United Kingdom},
Piotr Gwiazda
\footnote{Institute of Mathematics, Polish Academy of Sciences, 8 \'Sniadeckich Street, 00-656 Warsaw, Poland},
Karolina Kropielnicka
\footnote{Institute of Mathematics, Polish Academy of Sciences, 8 \'Sniadeckich Street, 00-656 Warsaw, Poland},
Anna Marciniak-Czochra
\footnote{Institute of Applied Mathematics, Interdisciplinary Center for scientific Computing (IWR) and BIOQUANT,  Heidelberg University, 69120 Heidelberg, Germany}
}   
   
\maketitle

\begin{abstract}
The Escalator Boxcar Train (EBT) method is a well known and widely used numerical method for one-dimensional structured population models of McKendrick--von Foerster type. Recently the method, in its full generality, has been applied to aged--structured two--sex population model (Fredrickson--Hoppensteadt model),  which consists of three coupled hyperbolic partial differential equations with nonlocal boundary conditions. We derive the {\it simplified} EBT method and prove its convergence to the solution of Fredrickson--Hoppensteadt model. The convergence can be proven, however only if we analyse the whole problem in the space of nonnegative Radon measures equipped with bounded Lipschitz distance (flat metric). Numerical simulations are presented to illustrate the results.
\end{abstract}

\section{Introduction}

Escalator Boxcar Train (EBT) method is a numerical integrator that was introduced in \cite{deRoos89} for a structured population models of McKendrick--von Foerster type \cite{McKendrick26} given by a scalar hyperbolic partial differential equation. It has been widely used because it approximates the density--function  of distribution of individuals in a way that has a clear biological interpretation: The method is based on representing the solution as a sum of masses localised in discrete points and tracing its dynamics along the characteristic lines of the model. 

The method has been popular for many years, even its convergence was proven only recently in \cite{EBT}, and in  \cite{GwJabMarUli14} where the rate of convergence was also shown. Results on the convergence were obtained using a theoretical approach to stability, where the underlying model is embedded in a space of nonnegative Radon measures $(\mathcal{M}^+(\bR_+))$ equipped with a bounded Lipschitz distance (flat metric). This approach was proposed in \cite{GLMC,GMC}. Such an {\em external approximation}, that is approximation of functions from ${\mathbf{L^1}}(\bR_+)$ by objects from a space of nonnegative Radon measures, is a natural consequence of the way in which the initial condition of the problem is approximated.

In this paper we focus on application of the EBT method to a system of structured population equations on example of the Fredrickson-Hoppensteadt model that is a  two--sex population model describing evolution of males and females and the process of heterogenous couples formation. The model was originally formulated in \cite{Fredrickson71} and later developed in \cite{Hoppensteadt75}. It consists of three population equations with structure, which are coupled through nonlocal boundary terms and a nonlocal and nonlinear source. Dynamics of males and females is given by McKendrick type equations, that is they consist of a transport equation with a growth term and  boundary terms determining the influx of newborn individuals. Evolution of couples is modelled by a similar equation, however it is equipped with so called marriage function. The marriage function describes the influx of new couples between males and females in the particular ages and at a certain moment. In reality, formation of new couples depends on many social and economical factors , such as religion, culture, education or health, thus it is a much more complicated process than birth and death rates for males and females. In the literature mentioned above, authors assume that the distribution of population is given by a density, so in the system (\ref{the_model_densities_nonlinear}), presented below, functions $u^m(t,x)$ and  $u^f(t,y)$ describe the distribution of males and females at time $t$ and age $x$ and $y$, respectively, while $u^c(t,x,y)$ is the number of couples at time $t$ between males at age $x$ and females at age $y$. The following system of nonlinear equations describes dynamics of the population of males, females and couples
\begin{eqnarray}
\partial_t u^m(t,x)+\partial_x u^m(t,x)+c^m(t,u^m(t,\cdot),u^f(t,\cdot),x)u^m(t,x)&=&0,\nonumber\\
u^m(t,0)&=&\int_{\bR_+^2}b^m(t,u^m(t,x),u^f(t,y),x,y)u^c(t,x,y)dxdy,\nonumber\\
u^m(0,x)&=& u^m_0(x),\nonumber\\\nonumber
\\\nonumber
\partial_t u^f(t,y)+\partial_y u^f(t,y)+c^f(t,u^m(t,\cdot),u^f(t,\cdot),y)u^f(t,y)&=&0,\nonumber\\
u^f(t,0)&=&\int_{\bR_+^2}b^f(t,u^m(t,x),u^f(t,y),x,y)u^c(t,x,y)dxdy,\nonumber\\
u^f(0,x)&=& u^f_0(x),\nonumber\\ \label{the_model_densities_nonlinear}
\end{eqnarray}
\begin{eqnarray*}
\partial_t u^c(t,x,y)+\partial_{x} u^c(t,x,y)+\partial_{y}u^c(t,x,y)+
c^c(t,u^m(t,\cdot),u^f(t,\cdot),u^
c(t,\cdot),x,y)u^c(t,x,y)&=&T(t,x,y),\\
u^c(t,x,0)=u^c(t,0,y)&=&0,\\
u^c(0,x,y)&=& u^c_0(x,y).\\
\end{eqnarray*}
Functions $c^m$, $c^f$ and $c^c$ describe the rates of disappearance of individuals, where disappearance of males or females is related to death, while couples disappearance reflects divorce or death of one of spouses. Functions $b^m$ and $b^f$ are birth rates of males and females. Observe that the mentioned coefficients  depend on ecological pressure in a nonlinear manner, that is they are nonlocal operators depending on the distribution of males, females and couples. 

The marriage function $T$ models the number of new marriages of males and females of age $x$ and $y$, respectively, at time $t$. It  also depends nonlinearly on the distribution of individuals. The choice of this function is a subject of ongoing discussions, see \cite{Hadeler89,Hadeler88,Martcheva99,Pruss94}, due to the properties like heterosexuality, homogeneity, consistency or competition. In this paper we follow the formulation proposed in \cite{Inaba93}, namely
\begin{align}
T(t,x,y)&=F(t,u^m(t,x),u^f(t,y),u^c(t,x,y),x,y)
\nonumber\\
&=\frac{\Theta(x,y)h(x)g(y)\left[u^m(t,x)-\int_0^\infty u^c(t,x,y)dy\right]\left[u^f(t,y)-\int_0^\infty u^c(t,x,y)dx\right]}{\gamma+\int_0^\infty h(x)\left[u^m(t,x)-\int_0^\infty u^c(t,x,y)dy\right]dx+\int_0^\infty g(y)\left[u^f(t,y)-\int_0^\infty u^c(t,x,y)dx\right]dy}.
\label{spec_marriage_function}
\end{align}

The function $\Theta(x,y)\in{\mathbf{L^1}}(\bR^2_+)\cap{\mathbf{L^\infty}}(\bR^2_+)$ describes the marriage rate of males of age $x$ and females of age $y$.  Notice that $\left[u^m(t,x)-\int_0^\infty u^c(t,x,y)dy\right]$ is the amount of unmarried males and $\left[u^f(t,y)-\int_0^\infty u^c(t,x,y)dx\right]$ is the number of unmarried females. The functions $h, g\in{\mathbf{L^1}}(\bR_+)\cap{\mathbf{L^\infty}}(\bR_+)$ describes the distribution of eligible males/females on the marriage market. We further assume that youngsters do not marry below a certain age $a$, i.e.
\begin{equation}\label{newbornsdonotmarry}
h(x)=g(y)=0 \mbox{ for } x,y\in [0,a_0)\,.
\end{equation}
The regularity of the remaining coefficients and their nonlinear dependences are presented in detail in Section \ref{definitions}.

The aim of the paper is analysis and convergence of a simplified EBT scheme, see Subsection 2.1, corresponding to (\ref{the_model_densities_nonlinear}). More precisely, we show convergence of the numerical integrator embedding the underlying problem (\ref{the_model_densities_nonlinear}) and its numerical scheme in a space of measures. The nonlinear age--structured, two--sex population model was presented and analysed in a space of nonnegative Radon measures in \cite{Ulikowska2012}, where the well possedness of the problem was proved.  

Recently, the EBT method has been derived in \cite{GKMC} for a partially linearised system of aged--structured equations of the Fredrickson-Hoppensteadt model. In such case, the only nonlinearity remains in function $T$ defined by (\ref{spec_marriage_function})):
\begin{eqnarray}
\partial_t u^m(t,x)+\partial_x u^m(t,x)+c^m(t,x)u^m(t,x)&=&0,\nonumber\\
u^m(t,0)&=&\int_{\bR_+^2}b^m(t,x,y)u^c(t,x,y)dxdy,\nonumber\\
u^m(0,x)&=& u^m_0(x),\nonumber\\\nonumber
\\\nonumber
\partial_t u^f(t,y)+\partial_y u^f(t,y)+c^f(t,y)u^f(t,y)&=&0,\nonumber\\
u^f(t,0)&=&\int_{\bR_+^2}b^f(t,x,y)u^c(t,x,y)dxdy,\nonumber\\
u^f(0,x)&=& u^f_0(x),\nonumber\\ \label{the_model_densities_linear}
\end{eqnarray}
\begin{eqnarray*}
\partial_t u^c(t,x,y)+\partial_{x} u^c(t,x,y)+\partial_{y}u^c(t,x,y)+c^c(t,x,y)u^c(t,x,y)&=&T(t,x,y),\\
u^c(t,x,0)=u^c(t,0,y)&=&0,\\
u^c(0,x,y)&=& u^c_0(x,y).\\
\end{eqnarray*}

Analysis of the partially linearised model  was an essential step towards the analysis of the nonlinear case (\ref{the_model_densities_nonlinear}), which we present in next sections.

The remainder of this paper is organised as follows. Section 2 is devoted to the EBT schemes. In Subsection 2.1 we present a simplified EBT method for the Fredrickson-Hoppensteadt model. In Subsection 2.2 we summarise the original EBT scheme for (\ref{the_model_densities_linear}) recently derived in \cite{GKMC}, while in Subsection 2.3 we show how to obtain a simplified EBT method from the original EBT approach \cite{GKMC}. The simplification of the method consists in unifying the rules (ODEs) in all cohorts, without distinction between internal and boundary ones. In Section 3, the underlying problem (\ref{the_model_densities_nonlinear}) is reformulated as an evolution in a space of nonnegative Radon measures. We approximate solutions of this problem  with a linear combination of Dirac Deltas, where the masses and localisations are obtained from the simplified EBT scheme embedded in the space of nonnegative Radon measures as well. Subsection \ref{definitions} is devoted to the analytical framework in the space of nonnegative Radon measures, where we introduce the necessary notation, definitions, lemmas and assumptions. The choice of state space allows proving the rate of convergence of the simplified EBT scheme in Section 4. Section 5 is devoted to numerical illustrations. As the computational error measurement is not trivial in flat metric (especially in $\bR^2$), some necessary details are provided in Subsection 5.1, while the obtained rate of convergence is illustrated in Subsections 5.2 and 5.3.


\section{Numerical methods based on EBT approach}
\setcounter{equation}{0}

In this paper we introduce and analyse a simplified EBT method for (\ref{the_model_densities_linear}). This Section starts with the presentation of the simplified method. In the later subsections we present the original method derived in \cite{GKMC}, and explain derivation of the simplified method from the original one.

\subsection{The simplified EBT method}

\begin{figure}[ht!]
	\centering
	\subfigure[Cohorts and internalisation moments for male population.]{
		\includegraphics[width=0.45\textwidth]{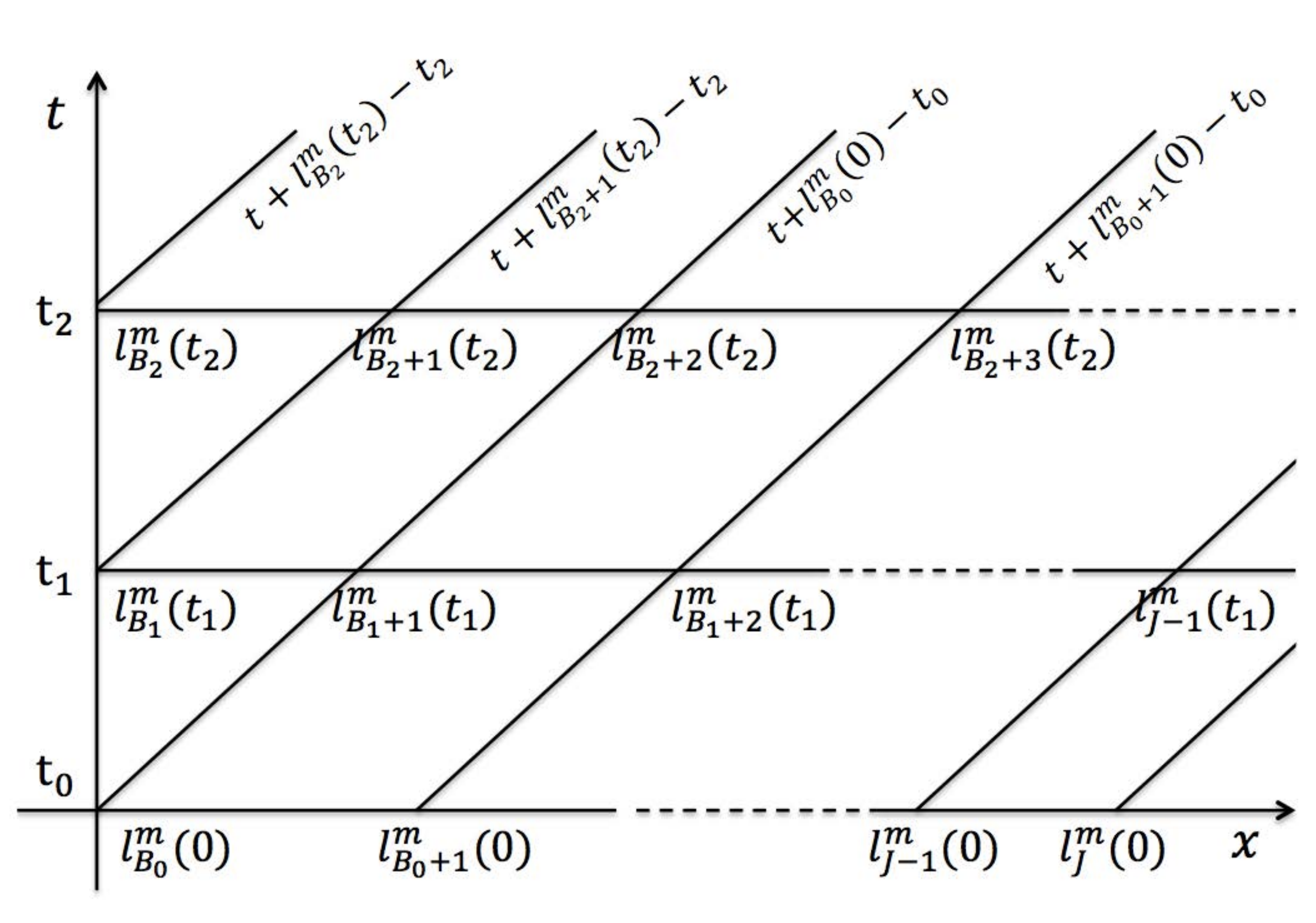}
		\label{cohorts_ini}
	}
	\hspace{1em}
	\subfigure[Masses and localisations of male population evolving in time interval $[t_n,t_{n+1})$.]{
		\includegraphics[width=0.45\textwidth]{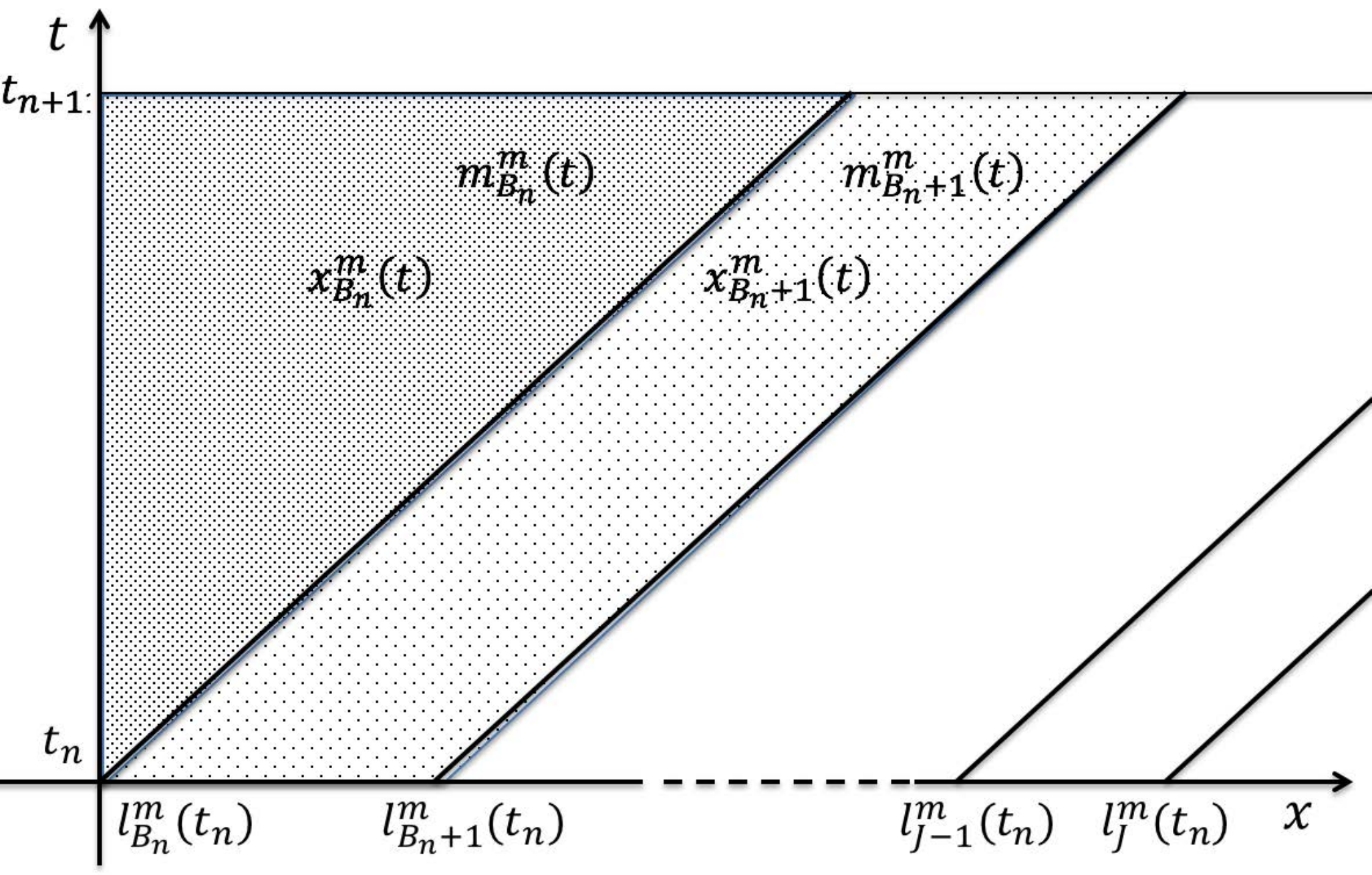}
	\label{cohorts}
	}
\end{figure}

As the concept of particle methods is grouping individuals into so called cohorts and tracing their dynamics in time, the first step of EBT algorithm consists in imposing, at time $t=0$, initial $J-B_0$ cohorts for males, females and 
$(J-B_0)^2$ for couples (see Figure \ref{cohorts_ini}) :
\begin{equation}\label{initial_internal_cohorts}
[l_i^m(0),l_{i+1}^m(0)),\ [l_j^f(0),l_{j+1}^f(0))\ {\rm and}\ [l_i^m(0),l_{i+1}^m(0))\times[l_j^f(0),l_{j+1}^f(0)),\ i,j=B_0,\ldots,J-1,
\end{equation}
respectively, in such a way, that ${\rm supp}(u_0^m)\subset [l_{B_0}^m(0),l_{J}^m(0))\ {\rm and}\ {\rm supp}(u_0^f)\subset [l_{B_0}^f(0),l_{J}^f(0)).$
Cohorts evolve in time along the characteristic lines of relevant transport operators in (\ref{the_model_densities_linear}). As always in a case of age--structured  problems, those characteristics are straight lines, see Figure \ref{cohorts_ini}.

In the next step we impose a mesh on the time variable $t\in[0,T)$ in such a way, that $t_0=0$ and $\bigcup_{n=0}^{N_T}[t_n,t_{n+1})=[0,T)$, see Figure \ref{cohorts_ini}. 

\begin{Remark}
In each time step $t_n$, new boundary cohorts are created. Boundary cohorts account for the influx of new males and females. In case of couples new boundary cohorts also appear, but they are empty, as we do not expect newborns to form couples due to \eqref{newbornsdonotmarry}. In Figure \ref{cohorts} we illustrate cohorts and internalisation moments for the male population during one time step. 
Notice that neither time steps nor boundaries of cohorts are to be equidistant.
\end{Remark}

\begin{Remark}
Depending on time $t$, the amount of cohorts changes. Let us notice, that in the internalisation moment $t_n$ male and female individuals are grouped into $J-B_n$ cohorts, while couples into $(J-B_n)^2$, that is
$$
[l_i^m(t_n),l_{i+1}^m(t_n)),\ i=B_n,\ldots,J,\ \hspace{1cm}
[l_j^f(t_n),l_{j+1}^f(t_n)),\ j=B_n,\ldots,J
$$
$$
[l_i^m(t_n),l_{i+1}^m(t_n))\times [l_j^f(t_n),l_{j+1}^f(t_n)),\ i,j=B_n,\ldots,J,
$$
If $t\in(t_k,t_{n+1})$ then we deal with $(J-B_n+1)$ cohorts for males and females and with $(J-B_n+1)^2$ cohorts for couples. Obviously individuals from time $t_n$ stay in their corresponding cohorts for time $t\in(t_n,t_{n+1})$:
$$
[t+l_i^m(t_n)-t_n,t+l_{i+1}^m(t_n)-t_n),\ i=B_n,\ldots,J,\ \hspace{1cm}
[t+l_j^f(t_n)-t_n,t+l_{j+1}^f(t_n)-t_n),\ j=B_n,\ldots,J
$$
$$
[t+l_i^m(t_n)-t_n,t+l_{i+1}^m(t_n)-t_n)\times [t+l_j^f(t_n)-t_n,t+l_{j+1}^f(t_n)-t_n),\ i,j=B_n,\ldots,J,
$$
but also new boundary cohorts appear (see the triangles in Figure \ref{cohorts}), that is
$$
[0,t+l_{B_n}^m(t_n)-t_n),\ \hspace{0.5cm}
[0,t+l_{B_n}^f(t_n)-t_n),\ \hspace{0.5cm}
[0,t+l_{B_n}^m(t_n)-t_n)\times [0,t+l_{B_n}^f(t_n)-t_n).
$$
Observe that the indexes of the cohorts at each of internalisation points move backward, i.e. $B_{n+1}=B_n-1$, with $B_0=0$ typically.
\end{Remark}

The idea of EBT algorithm is to trace the amount of individuals in cohorts, as the time changes from the moment $0$ to moment $T$, $t\in[0,T)$. Thus each cohort is characterised by  the mass and its location (see Figure \ref{cohorts}):
 $$
 (m_i^m(t),x_i^m(t)),\  (m_j^f(t),x_j^f(t))\ {\rm and}\ \left(m_{ij}^c(t),(x_{ij}^c(t),y_{ij}^c(t))\right).\ 
 $$
These quantities evolve in time intervals $[t_n,t_{n+1})$, $n=0,\ldots,{N_T}-1$ and are governed by the following system of differential equations, which constitutes the simplified EBT scheme, 
$$
\left\{
\begin{array}{rcl}\vspace{2mm}
\frac{d}{dt}m_i^m(t)&=&-c^m(t,x_i^m(t))m_i^m(t),\ i=B_n+1,\ldots,J\\ \vspace{2mm}
\frac{d}{dt}x_i^m(t)&=&1,\ i=B_n,\ldots,J\\ \vspace{2mm}
\frac{d}{dt}m_{B_n}^m(t)&=&-c^m(t,x_{B_n}^m(t))m_{B_n}^m(t)+\sum_{i,j={B_n}}^Jb^m(t,x_{ij}^c(t),y_{ij}^c(t))m_{ij}^c(t),\\
\end{array}
\right.
$$

\begin{equation}\label{the_numerical_method_simplified}
\left\{
\begin{array}{rcl}\vspace{2mm}
\frac{d}{dt}m_j^f(t)&=&-c^f(t,y_j^f(t))m_j^f(t),\ j=B_n+1,\ldots,J\\ \vspace{2mm}
\frac{d}{dt}y_j^f(t)&=&1,\ j=B_n,\ldots,J\\ \vspace{2mm}
\frac{d}{dt}m_{B_n}^f(t)&=&-c^f(t,x^f_{B_n}(t))m_{B_n}^f(t)+\sum_{i,j=B_n}^Jb^f(t,x_{ij}^c(t),y_{ij}^c(t))m_{ij}^c(t),
\end{array}
\right.
\end{equation}
$$
\left\{
\begin{array}{rcl}\vspace{2mm}
\frac{d}{dt}m_{ij}^c(t)&=&-c^c(t,x_{ij}^c(t),y_{ij}^c(t))m_{ij}^c(t)+\frac{N_{ij}(t)}{D_{ij}(t)}, \ i,j=B_n,\ldots,J\\ \vspace{2mm}
\frac{d}{dt}({x}_{ij}^c(t),{y}_{ij}^c(t))&=&(1,1), \ i,j=B_n,\ldots,J\\ 
\end{array}
\right.
$$
with
\begin{equation}\label{marriage_function_simplified}
\frac{N_{ij}(t)}{D_{ij}(t)}=
\frac{\Theta(x_{ij}^c(t),y_{ij}^c(t))h(x_{ij}^c(t))g(y_{ij}^c(t))\left(m_{i}^m(t)-\sum_{w=B_n}^Jm_{iw}^c(t)\right)\left(m_{j}^f(t)-\sum_{v=B_n}^Jm_{vj}^c(t)\right)}
{\gamma+\sum_{v=B_n}^Jh(x_{vj}^c(t))\left(m_{v}^m(t)-\sum_{w=B_n}^Jm_{vw}^c(t)\right)+
\sum_{w=B_n}^Jg(y_{iw}^c(t))\left(m_{w}^f(t)-\sum_{v=B_n}^Jm_{vw}^c(t)\right)}
\end{equation}

The definition of initial conditions for these equations depends on $n$:
\begin{enumerate}
\item[$(n=0)$]
Firstly let us observe that $t_0=0$, and the initial conditions should be consistent with the initial conditions of (\ref{the_model_densities_linear}). The boundary cohorts are defined with zero masses:
\begin{equation}\label{the_numerical_method_initial_1a}
(m_{B_0}^m(0),x_{B_0}^m(0))=(0,0),\  (m_{B_0}^f(0),y_{B_0}^f(0))=(0,0),
\end{equation}
\begin{equation}\label{the_numerical_method_initial_1b}
\left(m_{i j}^c(0),(x_{i j}^c(0),y_{i j}^c(0))\right)=(0,(x_i^m(0),y_j^f(0))),\   {\rm where}\ (i=B_0 \lor j=B_0) \land i,j\in\{B_0,\ldots,J\}.
\end{equation}
Initial conditions for the internal cohorts are derived from the biological model,
where masses mean the amount of individuals within a cohort:
\begin{equation}\label{the_numerical_method_initial_2}
\begin{split}
m_{i}^m(0)&=\int_{[l_{i-1}^m(0),l_{i}^m(0))}u_0^m(x)dx, \\
 m_{j}^f(0)&=\int_{[l_{j-1}^f(0),l_{j}^f(0))}u_0^f(y)dy,\\
m_{ij}^c(0)&=\int_{[l_{i-1}^m(0),l_{i}^m(0))\times[l_{j-1}^f(0),l_{j}^f(0))}u_0^c(x,y)dxdy,\quad i,j=B_0+1,\ldots,J,
\end{split}
\end{equation}
and locations mean the average value of structural variable within the underlying cohort:

\begin{equation}\label{the_numerical_method_initial_3}
\begin{split}
x_{i}^m(0)=&
\left\{
\begin{array}{l}
0\ {\rm if}\ m_i^m(0)=0,\\ 
\frac{1}{m_i^m(0)}\int_{[l_{i-1}^m(0),l_{i}^m(0))}xu_0^m(x)dx,\ {\rm otherwise},\\
\end{array}
\right.\\
y_{j}^f(0)=&
\left\{
\begin{array}{l}
0\ {\rm if}\ m_i^m(0)=0,\\ 
\frac{1}{m_j^f(0)}\int_{[l_{j-1}^f(0),l_{j}^f(0))}xu_0^f(y)dy,\ {\rm otherwise},\\
\end{array}
\right.\\
(x_{ij}^c(0),y_{ij}^c(0))=&
\quad (x_i^m(0),y_j^f(0)), \hspace{6cm} i,j=B_0+1,\ldots,J. \end{split}
\end{equation}

\item[$(n>0)$] Like previously, the boundary cohorts are defined with zero masses

\begin{equation}\label{the_numerical_method_initial_1a_k_step}
(m_{B_n}^m(t_n),x_{B_n}^m(t_n))=(0,0),\  (m_{B_n}^f(t_n),y_{B_n}^f(t_n))=(0,0),
\end{equation}
\begin{equation}\label{the_numerical_method_initial_1b_k_step}
\left(m_{i j}^c(t_n),(x_{i j}^c(t_n),y_{i j}^c(t_n))\right)=(0,(x_i^m(t_n),y_j^f(t_n))),\   {\rm where}\ (i=B_n \lor j=B_n) \land i,j\in\{B_n,\ldots,J\},
\end{equation}
while initial conditions for internal cohorts are obtained as an output of the $(n-1)$--th step of the algorithm (in the sense of limit $t\rightarrow t_{n}^-$):
\begin{equation}\label{the_numerical_method_initial_2_k_step}
\begin{split}
(x_i^m(t_n),m_i^m(t_n))=&\lim_{t\rightarrow t_{n}^-}(x_i^m(t),m_i^m(t)),\\
(y_j^f(t_n),m_j^f(t_n))=&\lim_{t\rightarrow t_{n}^-}(y_j^f(t),m_j^f(t)),\\
\left(m_{ij}^c(t_n),(x_{ij}^c(t_n),y_{ij}^c(t_n)))\right)=&\lim_{t\rightarrow t_{n}^-}\left(m_{ij}^c(t),(x_{ij}^c(t),y_{ij}^c(t))\right)
\end{split}
\end{equation}

\end{enumerate}

\subsection{The original EBT method}

The original EBT method for (\ref{the_model_densities_linear}) was derived in \cite{GKMC} and yields far more complicated form:
$$
\left\{
\begin{array}{rcl}\vspace{2mm}
\frac{d}{dt}m_i^m(t)&=&-c^m(t,x_i^m(t))m_i^m(t),\ i=B_n+1,\ldots,J\\ \vspace{2mm}
\frac{d}{dt}x_i^m(t)&=&1,\ i=B_n+1,\ldots,J\\ \vspace{2mm}
\frac{d}{dt}m_{B_n}^m(t)&=&-c^m(t,0)m_{B_n}^m(t)-\partial_xc^m(t,0)\Pi_{B_n}^m(t)\\
&&+\sum_{i,j={B_n}}^Jb^m(t,x_{ij}^c(t),y_{ij}^c(t))m_{ij}^c(t),\\[2mm]
\frac{d}{dt}\Pi_{B_n}^m(t)&=&m_{B_n}^m(t)-c^m(t,0)\Pi_{B_n}^m(t),\\[2mm]
x_{B_n}^m(t)&=&
\left\{
\begin{array}{l}
0\ {\rm if}\ m_{B_n}^m(t)=0,\\[2mm]
\frac{\Pi_{B_n}^m(t)}{m_{B_n}^m(t)},\ {\rm otherwise},\\
\end{array}
\right.\\
\end{array}
\right.
$$
\begin{equation}\label{the_numerical_method}
\left\{
\begin{array}{rcl}\vspace{2mm}
\frac{d}{dt}m_j^f(t)&=&-c^f(t,y_j^f(t))m_j^f(t),\ j=B_n+1,\ldots,J\\ \vspace{2mm}
\frac{d}{dt}y_j^f(t)&=&1,\ j=B_n+1,\ldots,J\\ \vspace{2mm}
\frac{d}{dt}m_{B_n}^f(t)&=&-c^f(t,0)m_{B_n}^f(t)-\partial_xc^f(t,0)\Pi_{B_n}^m(t)\\
&&+\sum_{i,j=B_n}^Jb^f(t,x_{ij}^c(t),y_{ij}^c(t))m_{ij}^c(t),\\[2mm]
\frac{d}{dt}\Pi_{B_n}^f(t)&=&m_{B_n}^f(t)-c^f(t,0)\Pi_{B_n}^f(t),\\[2mm]
y_{B_n}^f(t)&=&
\left\{
\begin{array}{l}
0\ {\rm if}\ m_{B_n}^f(t)=0,\\[2mm]
\frac{\Pi_{B_n}^m(t)}{m_{B_n}^f(t)},\ {\rm otherwise},\\
\end{array}
\right.\\
\end{array}
\right.
\end{equation}
$$
\left\{
\begin{array}{rcl}\vspace{2mm}
\frac{d}{dt}m_{ij}^c(t)&=&-c^c(t,x_{ij}^c(t),y_{ij}^c(t))m_{ij}^c(t)+\frac{N_{ij}(t)}{D_{ij}(t)}, \ i,j=B_n,\ldots,J\\[2mm]
\frac{d}{dt}(\tilde{x}_{ij}^c(t),\tilde{y}_{ij}^c(t))&=&\left[(1,1)-(x_{ij}^c(t),y_{ij}^c(y))c^c(t,x_{ij}^c(t),y_{ij}^c(t))\right]m_{ij}^c(t)+\frac{\bar{N}_{ij}(t)}{D_{ij}(t)},\\[2mm]
(x_{ij}^c(t),y_{ij}^c(t))&=&
\left\{
\begin{array}{l}
0\ {\rm if}\ m_{ij}^c(t)=0,\\[2mm]
\frac{(\tilde{x}_{ij}^c(t),\tilde{y}_{ij}^c(t))}{m_{ij}^c(t)},\ {\rm otherwise},\\
\end{array}
\right.\\
\end{array}
\right.
$$ 
where $\Pi_{B_n}^m(t)=\int_{{[l_{i-1}^m(t_n),l_{i}^m(t_n))}}xu^m(t,x)dx$ and $\Pi_{B_n}^f(t)=\int_{{[l_{j-1}^f(t_n),l_{j}^f(t_n))}}yu^f(t,y)dy$,\\
$$
(\tilde{x}_{ij}^c,\tilde{y}_{ij}^c)(t)=\int_{{[l_{i-1}^m(t_n),l_{i}^m(t_n))}\times {[l_{j-1}^f(t_n),l_{j}^f(t_n))}}(x,y)u^c(t,x,y)dxdy,
$$
\begin{eqnarray}\nonumber
N_{ij}(t)&=&\Theta(x_i^m(t),y_j^f(t))h(x_i^m(t))g(y_j^f(t))m_i^m(t)m_j^f(t)\\\label{Numerator}
&&-\sum_{v=B_n}^J\Theta(x_i^m(t),y_{vj}^c(t))h(x_i^m(t))g(y_{vj}^c(t))m_i^m(t)m_{vj}^c(t)\nonumber\\
&&-\sum_{w=B_n}^J\Theta(x_{iw}^c(t),y_j^f(t))h(x_{iw}^c(t))g(y_j^f(t))m_j^f(t)m_{iw}^c(t)\\ \nonumber\\
&&+\sum_{v,w=B_n}^J\Theta(x_{iw}^c(t),y_{vj}^c(t))h(x_{iw}^c(t))g(y_{vj}^c(t))m_{vj}^c(t)m_{iw}^c(t),\nonumber
\end{eqnarray}
\begin{eqnarray}\label{Numerator_bar}
\bar{N}_{ij}(t)&=&(x_i^m(t),y_j^f(t))\Theta(x_i^m(t),y_j^f(t))h(x_i^m(t))g(y_j^f(t))m_i^m(t)m_j^f(t)\\\nonumber
&&-\sum_{v=B_n}^J(x_i^m(t),y_{vj}^c(t)) \Theta(x_i^m(t),y_{vj}^c(t))h(x_i^m(t))g(y_{vj}^c(t))m_i^m(t)m_{vj}^c(t)\\\nonumber
&&-\sum_{w=B_n}^J(x_{iw}^c(t),y_j^f(t))(\Theta(x_{iw}^c(t),y_j^f(t))h(x_{iw}^c(t))g(y_j^f(t))m_j^f(t)m_{iw}^c(t)\\ \nonumber
&&+\sum_{v,w=B_n}^J(x_{iw}^c(t),y_{vj}^c(t))\Theta(x_{iw}^c(t),y_{vj}^c(t))h(x_{iw}^c(t))g(y_{vj}^c(t))m_{vj}^c(t)m_{iw}^c(t)\nonumber
\end{eqnarray}
and
\begin{eqnarray}\label{Denominator}
D_{ij}(t)&=&\gamma+\sum_{i=B_n}^Jh(x_i^m(t))m_i^m(t)-\sum_{i,j=B_n}^Jh(x_{ij}^c(t))m_{ij}^c(t)\\ \nonumber
&&+\sum_{j=B_n}^Jg(y_j^f(t))m_j^f(t)-\sum_{i,j=B_n}^Jg(y_{ij}^c(t))m_{ij}^c(t).
\end{eqnarray}

The functions $N_{ij}(t),\ \bar{N}_{ij}(y)$ and $D_{ij}(t)$ appear in the EBT scheme, due to the generic marriage function (\ref{spec_marriage_function}). Namely, the integral of $F$ over a cohort is approximated as
$$
\int_{[l_{i-1}^m(t_k),l_{i}^m(t_k))\times[l_{j-1}^f(t_k),l_{j}^f(t_k))}F(t,u^m(t,x),u^f(t,x),u^c(t,x,y),x,y)dxdy \simeq \frac{N_{ij}(t)}{D_{ij}(t)},
$$
while its first moment is approximated as
$$
\int_{[l_{i-1}^m(t_k),l_{i}^m(t_k))\times[l_{j-1}^f(t_k),l_{j}^f(t_k))}(x,y)F(t,u^m(t,x),u^f(t,x),u^c(t,x,y),x,y)dxdy \simeq \frac{\bar{N}_{ij}(t)}{D_{ij}(t)}.
$$ 
For further details we refer to \cite{GKMC}. Differential equations constituting the EBT method (\ref{the_numerical_method}) are equipped with the with initial boundary conditions (\ref{the_numerical_method_initial_1a})--(\ref{the_numerical_method_initial_2_k_step}), where instead of (\ref{the_numerical_method_initial_1b}) and (\ref{the_numerical_method_initial_1b_k_step}), the following constraints are applied
\begin{equation}\label{the_numerical_method_initial_1b_simplified_k_step}
\left(m_{i j}^c(t_k),(x_{i j}^c(t_k),y_{i j}^c(t_k))\right)=(0,(0,0)),\   {\rm where}\ (i=B_n \lor j=B_n) \land i,j\in\{B_n,\ldots,J\},
\end{equation}
for all $n$.

\subsection{Relation between the methods}

The simplified method has a far more transparent representation that the original EBT model for couples in \cite{GKMC}, it requires less computations while formally remaining the same order of convergence as we will show in the next section. In case of male and female populations it only differs by the evolution of localisations in the boundary cohorts. This simplification was already proposed in \cite{EBT} for the single species case, whose convergence was proven in \cite{GwJabMarUli14}. Here compared to \cite{GKMC}, we simplified the EBT scheme further by choosing the localisation of the couples in the boundary cohorts consistently in terms of the localisations of the male and female populations. More precisely, \eqref{the_numerical_method_initial_1b_simplified_k_step} is imposed compared to \eqref{the_numerical_method_initial_1b}-\eqref{the_numerical_method_initial_1b_k_step}. This results in simpler approximation of the marriage function and in the fact that the characteristics in age variable of couples, male and female populations remain the same.

In fact, we can check that under this consistent choice of the initial data for couples, female and male populations, the simplified EBT scheme is a particular case of the original EBT scheme. Notice, that the formula for $(\tilde{x}_{ij}^c(t),\tilde{y}_{ij}^c(t))$ in (\ref{the_numerical_method}) together with $(x_{ij}^c(t),y_{ij}^c(t))=\frac{1}{m_{ij}^c(t)}(\tilde{x}_{ij}^c(t),\tilde{y}_{ij}^c(t))$ yields
\begin{equation}\label{couples_boundary_con_1}
\frac{d}{dt}(x_{ij}^c,y_{ij}^c)(t)=\frac{d}{dt}\frac{(\tilde{x}_{ij}^c,\tilde{y}_{ij}^c)(t)}{m^c_{ij}(t)}=(1,1)-\frac{(x_{ij}^c(t),y_{ij}^c(t))}{m_{ij}^c(t)}\frac{N_{ij}(t)}{D_{ij}(t)}+\frac{\bar{N}_{ij}(t)}{m_{ij}^c(t)D_{ij}(t)}
\end{equation}
for all $t\in [t_n,t_{n+1})$. Let us remind that $\frac{d}{dt}x_i^m(t)=1$ and analyse only the first component of the equality (\ref{couples_boundary_con_1}), as the same reasoning can be made for $y_{ij}^f(t)$. We can rewrite this ODE as 
\begin{align*}
\frac{d}{dt}\left(x_{ij}^c(t)-x_i^m(t)\right)
=&\frac{1}{ m_{ij}^c(t) D_{ij}(t)}\left[{\bar{N}}_{ij}(t)\cdot \left(
\begin{array}{c}\vspace{1.5mm}
1\\
0\\
\end{array}
\right)
- x_{ij}^c(t) N_{ij}(t)\right]
\\
=&\frac{1}{ m_{ij}^c(t) D_{ij}(t)}
\Bigg[( x_i^m(t)- x_{ij}^c(t))\Theta( x_i^m(t), y_j^f(t))h( x_i^m(t))g( y_j^f(t)) m_i^m(t) m_j^f(t)\\
&-( x_i^m(t)- x_{ij}^c(t))\sum_{v=B_n}^J \Theta( x_i^m(t), y_{vj}^c(t))h( x_i^m(t))g( y_{vj}^c(t)) m_i^m(t) m_{vj}^c(t)\\\nonumber
&-\sum_{w=B}^J( x_{iw}^c(t)- x_{ij}^c(t))(\Theta( x_{iw}^c(t), y_j^f(t))h( x_{iw}^c(t))g( y_j^f(t)) m_j^f(t) m_{iw}^c(t)\\ \nonumber
&+\sum_{v,w=B}^J( x_{iw}^c(t)- x_{ij}^c(t))\Theta( x_{iw}^c(t), y_{vj}^c(t))h( x_{iw}^c(t))g( y_{vj}^c(t)) m_{vj}^c(t) m_{iw}^c(t)\Bigg]
\end{align*}
\begin{align*}
\hspace{4.5cm}=&\frac{1}{ m_{ij}^c(t) D_{ij}(t)}
\Bigg[( x_i^m(t)- x_{ij}^c(t))\Theta( x_i^m(t), y_j^f(t))h( x_i^m(t))g( y_j^f(t)) m_i^m(t) m_j^f(t)\\
&-( x_i^m(t)- x_{ij}^c(t))\sum_{v=B_n}^J \Theta( x_i^m(t), y_{vj}^c(t))h( x_i^m(t))g( y_{vj}^c(t)) m_i^m(t) m_{vj}^c(t)\\\nonumber
&-\sum_{w=B}^J( x_{iw}^c(t) \mp x_i^m(t) - x_{ij}^c(t))(\Theta( x_{iw}^c(t), y_j^f(t))h( x_{iw}^c(t))g( y_j^f(t)) m_j^f(t) m_{iw}^c(t)\\ \nonumber
&+\sum_{v,w=B}^J( x_{iw}^c(t) \mp x_i^m(t) - x_{ij}^c(t))\Theta( x_{iw}^c(t), y_{vj}^c(t))h( x_{iw}^c(t))g( y_{vj}^c(t)) m_{vj}^c(t) m_{iw}^c(t)\Bigg].
\end{align*}
One can rewrite the previous system in terms of the auxiliary variables $z_{ij}(t)=x_{ij}^c(t)-x_i^m(t)$. Let us assume that $x_i^m(t_n)=x_{ij}^c(t_n)$ or equivalently $z_{ij}(t_n)=0$ for $j=B_n,\ldots,J$. It is easy to observe from the previous expression that $z_{ij}(t)=0$ is a solution of the system of ODEs consistent with the initial conditions at $t_n$, and thus by the uniqueness of the ODE system we deduce that $z_{ij}(t)=x_{ij}^c(t)-x_i^m(t)=0$ for all $t\in [t_n,t_{n+1})$. As a consequence, we infer that 
$$
\frac{d}{dt}\left(x_{ij}^c(t)-t\right)=0,
$$
which explains the last formula in the simplified EBT ODE system (\ref{the_numerical_method_simplified}). Notice, that in this case, when $x_i^m(t)=x_{ij}^c(t)$ and $y_j^f(t)=y_{ij}^c(t)$ for $i,j=B(t),\ldots,J$, expression $\frac{N_{ij}(t)}{D_{ij}(t)}$ simplifies significantly as \eqref{marriage_function_simplified}.


\section{Embedding in a space of measures}
\setcounter{equation}{0}

As it was stated in the introduction, we are going to describe and analyse underlying problems and their solutions in a space of nonnegative Radon measures equipped with flat metric. In fact, setting some models of population dynamics in this space was suggested for the first time in \cite{GLMC}. The aged--structured two-sex population model with age as a structure variable  was, actually, embedded in a suitable space in \cite{Ulikowska2012}, where the approach followed after \cite{CCGU12,GLMC,GMC}.

Alternalively one can investigate the underlying problem not in the positive cone in the space of measures 
$\mathcal{M}_+(\bR_+^N)$ but in the full Banach space which is a closure of space of bounded Radon measures $\mathcal{M}(\bR_+^N)$ with respect to bounded Lipschitz distance, \cite{HilleWorm09}. See also \cite{Weaver99} for similar {\it Lipschitz--free space}. It is important to point out that the space defined in \cite{HilleWorm09} is predual to $W^{1,\infty}(\bR^N_+)$ and is essentially smaller then $\left(W^{1,\infty}(\bR^N_+)\right)^*$ .

The crucial reason to consider the predual space instead of $\left(W^{1,\infty}(\bR^N_+)\right)^*$ is the lack of the continuity of the semigroup generated by the transport operator, see Lemma 2 in \cite{GOU}.

So through this paper,   $\mathcal{M}^+(\bR^i_+)$ denotes the space of nonnegative Radon measures  with bounded total variation, where $\bR^i_+=[0,\infty)^i,\ i=1,2$, and $B\in\mathcal{B}(\bR_+)$ is a Borel set. We will investigate the following equivalent of system (\ref{the_model_densities_nonlinear}):

\begin{equation}\label{the_model_measures}
\begin{array}{rcl}
\partial_t\mu_t^m+\partial_x\mu_t^m+\xi^m(t,\mu_t^m,\mu_t^f)\mu_t^m&=&0,\hspace{1.4cm}(t,x)\in[0,T]\times \bR_+\\
D_\lambda\mu_t^m(0^+)&=&\int_{\bR_+^2}\beta^m(t,\mu_t^m,\mu_t^f)(z)d\mu_t^c(z)\\
\mu_0^m& \in & \mathcal{M}^+(\bR_+)\\
\\
\partial_t\mu_t^f+\partial_x\mu_t^f+\xi^f(t,\mu_t^m,\mu_t^f)\mu_t^f&=&0,\hspace{1.4cm}(t,x)\in[0,T]\times \bR_+\\
D_\lambda\mu_t^f(0^+)&=&\int_{\bR_+^2}\beta^f(t,\mu_t^m,\mu_t^f)(z)d\mu_t^c(z)\\
\mu_0^f& \in & \mathcal{M}^+(\bR_+)\\
\\
\partial_t\mu_t^c+\partial_{z_1}\mu_t^c+\partial_{z_2}\mu_t^c+\xi^c(t,\mu_t^m,\mu_t^f,\mu_t^c)\mu_t^c&=&\mathcal{T}(t,\mu_t^m,\mu_t^f,\mu_t^c),\hspace{4mm}(t,z)\in[0,T]\times \bR_+^2\\
\mu_t^c(\{0\}\times B)=\mu_t^c(B\times \{0\})&=&0\\
\mu_0^m& \in & \mathcal{M}^+(\bR^2_+)\\

\end{array}
\end{equation}

Measures $\mu_t^m$, $\mu_t^f\in \mathcal{M}_+(\bR^1_+)$ and $\mu_t^c\ \in \mathcal{M}_+(\bR^2_+)$ describe the distribution of males, females and couples, respectively, at time $t$. Functions $\xi^m$, $\xi^f$ and $\xi^c$ are equivalents of functions $c^m$, $c^f$ and $c^c$ and describe the disappearance of individuals, while functions $\beta^m$ and $\beta^f$ constitute equivalents for $b^m$ and $b^f$, birth rates of males and females, respectively. Symbols $D_\lambda\mu_t^m(0^+)$ and $D_\lambda\mu_t^f(0^+)$ denote Radon--Nikodym derivatives of $\mu_t^m$ and $\mu_t^f$, respectively, with respect to the one dimensional Lebesgue measure $\lambda$ at point $0$, as we assume that the support of singular part of measures $\mu_t^m$ and $\mu_t^f$ does not contain $0$. 

Before we make a comment on the marriage function, let us define the distribution of single males and females, $s_t^m$ and $s_t^f$, respectively. By {\it single} $s_t^m$ or $s_t^f$we mean not only those who has never been married by time $t$, but also those who are divorced or widowed at time $t$. Let measures $\sigma_t^m$ and $\sigma_t^f$ be projections of $\mu_t^c$ on $\bR_+$ and describe a distribution of males and females respectively, who are married at time $t$:
\begin{equation}\label{newlyweds}
\sigma_t^m(B)=\mu_t^c(B\times\bR_+),\ {\rm and}\ \sigma_t^f(B)=\mu_t^c(\bR_+\times B).
\end{equation}
Now $s_t^m$ and $s_t^f$ can be easily defined as
\begin{equation}\label{singles}
s_t^m(B)=(\mu_t^m-\sigma_t^m)(B\times\bR_+),\ {\rm and}\ s_t^f(B)=(\mu_t^f-\sigma_t^f)(B\times\bR_+)
\end{equation}

Following \cite{Ulikowska2012}  and \cite{GKMC} we adopt the following definition of generic marriage function:
\begin{equation}\label{spec_marriage_function_measures}
\mathcal{T}(t,\mu_t^m,\mu_t^f,\mu_t^c))=\mathcal{F}(t,\mu_t^m-\sigma_t^m,\mu_t^f-\sigma_t^f))=
\end{equation}
$$
\mathcal{F}(t,s_t^m,s_t^f)=\frac{\Theta(x,y)h(x)g(y)}{\gamma+\int_0^\infty h(z)ds_t^m(z)+\int_0^\infty g(w)ds_t^f(w)}(s_t^m\otimes s_t^f),
$$
where $(s_t^m\otimes s_t^f)$ is a product measure on $\bR_+^2$.

\subsection{EBT schemes for fully nonlinear model in space of measures}\label{EBT_measures}

The output of the numerical method should evolve in the same space that the solution of (\ref{the_model_measures}). This requirement can be easily satisfied by defining measures $\nu_{k,t}^m$, $\nu_{k,t}^f$ and $\nu_{k,t}^c$ as a linear combination of Dirac measures:

\begin{equation}\label{solution_of_numerical_method}
\nu_{k,t}^m:=\sum_{i=B_n}^Jm_i^m(t)\delta_{\{x_i^m(t)\}},\ \nu_{k,t}^f:=\sum_{j=B_n}^Jm_j^f(t)\delta_{\{y_j^f(t)\}},\ \nu_{k,t}^c:=\sum_{i,j=B_n}^Jm_{ij}^c(t)\delta_{\{x_{ij}^c(t),y_{ij}^c(t)\}},
\end{equation}
where $t\in[t_n,t_{n+1})$, (such that $t_{n+1}-t_n\leq a_0$), and $(x_i^m(t), m_i^m(t))$, $(y_j^f(t), m_j^f(t))$, $((x_{ij}^c(t),y_{ij}^c(t)), m_{ij}^c(t))$ is the output of the EBT algorithm (\ref{the_numerical_method}). The subscript $k$ ($k=J-B_0$)  is related to the approximation of initial condition, as it is equal to the amount of initial cohorts for males and females,
\begin{equation}\label{initial_approximation}
\nu_{k,0}^m:=\sum_{i=B_0}^Jm_i^m(0)\delta_{\{x_i^m(0)\}},\ \nu_{k,0}^f:=\sum_{j=B_0}^Jm_j^f(0)\delta_{\{y_j^f(0)\}},\ 
\nu_{k,0}^c:=\sum_{i,j=B_0}^Jm_{ij}^c(0)\delta_{\{x_{ij}^c(0),y_{ij}^c(0)\}}.
\end{equation}
For the transparency of notation we will omit subscript $k$ in the further part of the manuscript.

For $t\in[t_n,t_{n+1})$ the system (\ref{the_model_measures}) is approximated with the following EBT scheme 
$$
\left\{
\begin{array}{rcl}\vspace{2mm}
\frac{d}{dt}m_i^m(t)&=&-\xi^m(t,\nu_t^m,\nu_t^f)(x_i^m(t))m_i^m(t),\ i=B_n+1,\ldots,J\\ \vspace{2mm}
\frac{d}{dt}x_i^m(t)&=&1,\ i=B_n,\ldots,J\\ \vspace{2mm}
\frac{d}{dt}m_{B_n}^m(t)&=&-\xi^m(t,\nu_t^m,\nu_t^f)(x_{B_n}^m(t))m_{B_n}^m(t)+\sum_{i,j={B_n}}^J\beta^m(t,\nu_t^m,\nu_t^f)(x_{ij}^c(t),y_{ij}^c(t))m_{ij}^c(t),\\
\end{array}
\right.
$$
\begin{equation}\label{the_numerical_method_measures}
\left\{
\begin{array}{rcl}\vspace{2mm}
\frac{d}{dt}m_j^f(t)&=&-\xi^f(t,\nu_t^m,\nu_t^f)(y_j^f(t))m_j^f(t),\ j=B_n+1,\ldots,J\\ \vspace{2mm}
\frac{d}{dt}y_j^f(t)&=&1,\ j=B_n,\ldots,J\\ \vspace{2mm}
\frac{d}{dt}m_{B_n}^f(t)&=&-\xi^f(t,\nu_t^m,\nu_t^f)(y_{B_n}^f(t))m_{B_n}^f(t)+\sum_{i,j=B_n}^J\beta^f(t,\nu_t^m,\nu_t^f)(x_{ij}^c(t),y_{ij}^c(t))m_{ij}^c(t),\\
\end{array}
\right.
\end{equation}

$$
\left\{
\begin{array}{rcl}\vspace{2mm}
\frac{d}{dt}m_{ij}^c(t)&=&-\xi^c(t,\nu_t^m,\nu_t^f,\nu_t^c)(x_{ij}^c(t),y_{ij}^c(t))m_{ij}^c(t)+\frac{N_{ij}(t)}{D_{ij}(t)}, \ i,j=B_n,\ldots,J\\ \vspace{2mm}
\frac{d}{dt}({x}_{ij}^c(t),{y}_{ij}^c(t))&=&(1,1),\ i,j=B_n,\ldots,J\\
\end{array}
\right.
$$
where $\frac{N_{ij}(t)}{D_{ij}(t)}$ is defined by (\ref{marriage_function_simplified}).

\begin{Remark}\label{remark_marriage_by_dirac}
Let us notice that  application of  measures (\ref{solution_of_numerical_method}) to function $\mathcal{T}$ (defined in (\ref{spec_marriage_function_measures})) results in the formula:
$$
\mathcal{T}(t,\nu_t^m,\nu_t^f,\nu_t^c)(x_{ij}^c(t),y_{ij}^c(t))=\sum_{i,j=B_n}^J\frac{N_{ij}(t)}{D_{ij}(t)}\delta_{\{x_{ij}^c(t),y_{ij}^c(t)\}},
$$ 
for $t\in[t_n,t_{n+1})$.
\end{Remark}

\subsection{Notation, definitions and important facts}\label{definitions}

Through this paper we understand, that $\mathbf{C}^1(\bR_+^N;\bR)$ is a class of differentiable functions from $\bR_+^N$ to $\bR$.
Defined below Lipschitz bounded distance (flat metric) $d_N$ is associated, in the paper, with space $\mathcal{M}_+(\bR_+^N)$, where $N\in\bN$.
\begin{Definition}\label{the_flat_metric} 
Let $\mu, \nu \in{\mathcal M_+}({\bR^N_{+}})$, where $N\in\bN$.  The distance
function $d_N : {\mathcal M_+}({\bR^N_{+}}) \times  {\mathcal M_+}({\bR^N_{+}}) \rightarrow
[0, \infty)$ is defined by
\begin{equation*}
d_N(\mu_1,\mu_2)= \sup\left\{\int_{\bR_+^N}\varphi\,d(\mu_1-\mu_2):\ \varphi\in\mathbf{C}^1(\bR_+^N;\bR)\ {\rm and}\ \|\varphi\|_{\mathbf{W}^{\mathbf{1},\infty}(\bR_+^N,\bR)}\leq 1\right\},
\end{equation*}
where $\|\varphi\|_{\mathbf{W}^{\mathbf{1},\infty}(\bR_+^N,\bR)}=\max\{\|\varphi\|_{\mathbf{L}^\infty},\|\partial_{x_1}\varphi\|_{\mathbf{L}^\infty},\|\partial_{x_2}\varphi\|_{\mathbf{L}^\infty},\ldots,\|\partial_{x_N}\varphi\|_{\mathbf{L}^\infty}\}$, $N\in\bN$.
\end{Definition}

The metric $d_N$ is a distance derived from the dual norm of $\mathbf{W}^{\mathbf{1},\infty}(\bR_+^N,\bR)$. For sake of simplicity we will write $\|\cdot\|_{\mathbf{W}^{\mathbf{1},\infty}}$ instead of 
$\|\cdot\|_{\mathbf{W}^{\mathbf{1},\infty}(\bR_+^N,\bR)}$ unless it leads to misunderstandings. Space $\left( \mathcal{M}_+(\bR^N_+),d_n\right)$ possesses two favourable features which are necessary for investigating the convergence of numerical scheme. It is well-known, see \cite{Spohn,CCC,GMC}, that the metric space $\left( \mathcal{M}_+(\bR^N_+),d_N\right)$ is complete separable and that the convergence in metric $d_N$ is equivalent to narrow convergence of sequences in $\mathcal{M}_+(\bR^N_+)$. 

As we are concerned with linear combinations of Dirac measures it is worthy to make the following observation.

\begin{Lemma}\label{LemmaFlatMetric}
Let $\mu = \sum_{i=1}^{J}m_i\delta_{x_i}$ and  $\tilde \mu = \sum_{i=1}^{J} \tilde m_i \delta_{\tilde x_i}$, where $J\in\bN$, $x_i$, $\tilde x_i\in\bR_+^N$ and $m_i$, $\tilde m_i\in\bR_+$. Then,
\begin{eqnarray*}
d_N\left(\mu, \tilde \mu \right)
&\leq&
\sum_{i=1}^{J} \left( \|x_i - \tilde x_i\|m_i + |m_i - \tilde m_i| \right).
\end{eqnarray*}
\end{Lemma}

\begin{proofof}
The proof consists in applying the triangle inequality and Definition  \ref{the_flat_metric} to get
\begin{eqnarray*}
d_N\left(\mu, \tilde \mu \right)
&\leq&
d_N\left( \sum_{i=1}^{J} m_i \delta_{x_{i}},  \sum_{i=1}^{J} m_i \delta_{\tilde x_{i}} \right)
+
d_N\left( \sum_{i=1}^{J} m_i \delta_{\tilde x_{i}},  \sum_{i=1}^{J} \tilde m_i \delta_{\tilde x_{i}}\right)\\
&\leq&
\sup\left\{\int_{\bR_+^N}\varphi \,d\left( \sum_{i=1}^{J} m_i \delta_{x_{i}}-  \sum_{i=1}^{J} m_i \delta_{\tilde x_{i}} \right):\ \varphi\in\mathbf{C}^1(\bR_+^N;\bR)\ {\rm and}\|\varphi\|_{\mathbf{W}^{\mathbf{1},\infty}(\bR_+^N,\bR)}\leq 1\right\}\\
&&+
\sup\left\{\int_{\bR_+^N}\varphi \,d\left( \sum_{i=1}^{J} m_i \delta_{\tilde x_{i}}-  \sum_{i=1}^{J} \tilde m_i \delta_{\tilde x_{i}}\right):\ \varphi\in\mathbf{C}^1(\bR_+^N;\bR)\ {\rm and}\|\varphi\|_{\mathbf{W}^{\mathbf{1},\infty}(\bR_+^N,\bR)}\leq 1\right\}\\
&\leq&
\sup\left\{ \sum_{i=1}^{J} |\varphi(x_i)-\varphi(\tilde x_i)|m_i :\ \varphi\in\mathbf{C}^1(\bR_+^N;\bR)\ {\rm and}\|\varphi\|_{\mathbf{W}^{\mathbf{1},\infty}(\bR_+^N,\bR)}\leq 1\right\}\\
&&+
\sup\left\{\sum_{i=1}^{J} \varphi(\tilde x_{i}) |m_i-\tilde m_i| :\ \varphi\in\mathbf{C}^1(\bR_+^N;\bR)\ {\rm and}\|\varphi\|_{\mathbf{W}^{\mathbf{1},\infty}(\bR_+^N,\bR)}\leq 1\right\}\\
&\leq&
\sum_{i=1}^{J} \|x_i-\tilde x_i\|m_i +\sum_{i=1}^{J}|m_i-\tilde m_i|.
\end{eqnarray*}
\end{proofof}

\begin{Definition}\label{semiflow}
Let $(E,d)$ be a metric space. A family of bounded operators $S:[0,T]\times[0,T]\times E\rightarrow E$ is called a Lipschitz semiflow if for steps $s,t\in[0,T]$, and times $\tau$ such that $\tau,\tau+s,\tau+t,\tau+s+t\in[0,T]$ the following conditions are satisfied,
\begin{enumerate}
\item
$S(0;\tau)=I$,
\item
$S(t+s;\tau)=S(t;\tau+s)S(s;\tau)$,
\item
$d(S(t;\tau)\mu,S(s;\tau)\nu)\leq L(d(\mu,\nu)+|t-s|)$,
\end{enumerate}
\end{Definition}
\smallskip
Let us define product spaces
$$
\mathcal{U}=\mathcal{M}_+(\bR_+)\times \mathcal{M}_+(\bR_+)\times \mathcal{M}_+(\bR_+^2)\ {\rm and}\ \mathcal{V}=\mathcal{M}_+(\bR_+)\times \mathcal{M}_+(\bR_+).
$$
We will investigate the problem of approximation of weak solutions of model (\ref{the_model_measures}) in the metric space $(\mathcal{U},\mathbf{d})$, where $\mathbf{d}=d_1+d_1+d_2$.

The following definition of weak solution to (\ref{the_model_measures}) was proposed in \cite{Ulikowska2012} (see also \cite{CCGU12}).

\begin{Definition}\label{weak_solution}
A triple $\mathbf{u}=(\mu^m,\mu^f,\mu^c):[a,b]\rightarrow \mathcal{U}$ is a weak solution to the system (\ref{the_model_measures}) on the interval $[a,b]$, if $\mu^m,\mu^f,\mu^c$ are narrowly continuous with respect to time and for all $(\varphi^m,\varphi^f,\varphi^c)$ such that $\varphi^m,\varphi^f\in(\mathbf{C^1}\cap\mathbf{W^{1,\infty}})([a,b]\times\bR_+;\bR)$ and $\varphi^c\in(\mathbf{C^1}\cap\mathbf{W^{1,\infty}})([a,b]\times\bR_+^2;\bR)$, the following equalities hold
$$
\int_a^b\int_{\bR_+}\left(\partial_t\varphi^i(t,x)+\partial_x\varphi^i(t,x)-\xi^i(t,\mu_t^m,\mu_t^f)\varphi^i(t,x)\right)d\mu_t^i(x)dt
+\int_a^b\varphi^i(t,0)\int_{\bR_+^2}\beta^i(t,\mu_t^m,\mu_t^f)(z)d\mu_t^c(z)dt
$$
$$
=\int_{\bR_+}\varphi^i(b,x)d\mu_b^i(x)-\int_{\bR_+}\varphi^i(a,x)d\mu_a^i(x),\ {\rm for}\ i=f,m\
$$
and
$$
\int_a^b\!\!\int_{\bR_+}\!\!\left(\partial_t\varphi^c(t,z)+\partial_x\varphi^c(t,z)+\partial_y\varphi^c(t,z)-\xi^c(t,\mu_t^m,\mu_t^f,\mu_t^c)\varphi^c(t,z)\right)d\mu_t^c(z)dt
+\int_a^b\!\!\int_{\bR_+^2}\!\!\varphi^c(t,z)d\mathcal{T}(t,\mu_t^m,\mu_t^f,\mu_t^c)(z)dt
$$
$$
=\int_{\bR_+}\varphi^c(b,z)d\mu_b^c(z)-\int_{\bR_+}\varphi^c(a,z)d\mu_a^c(z).
$$
\end{Definition}
Here, narrowly continuous functions are understood in a sense of narrow convergence introduced in \cite{AGS}. 

\begin{Assumption}\label{model_functions_measures}
We make the following assumptions on the model functions:
$$
\begin{array}{rcl}
\xi^m,\xi^f&\in&\mathbf{BC}^{0,1}([0,T]\times\mathcal{V},\mathbf{W^{1,\infty}}(\bR_+;\bR))\\
\beta^m,\beta^f&\in&\mathbf{BC}^{0,1}([0,T]\times\mathcal{V},\mathbf{W^{1,\infty}}(\bR_+^2;\bR))\\
\xi^c&\in&\mathbf{BC}^{0,1}([0,T]\times\mathcal{U},\mathbf{W^{1,\infty}}(\bR_+^2;\bR))\\
\mathcal{T}&\in&\mathbf{BC}^{0,1}([0,T]\times\mathcal{U},\mathcal{M}_+(\bR_+^2))\\
\end{array}
$$
We understand, that spaces $\mathbf{BC}^{0,1}([0,T]\times\mathcal{V},X)$ and $\mathbf{BC}^{0,1}([0,T]\times\mathcal{U},X)$ are spaces of $X$ valued functions, bounded with respect to the $\|\cdot\|_X$ norm, continuous with respect to time and Lipschitz continuous with respect to measure variables. We understand, that function $\mathcal{T}$ is bounded with respect to $\|\cdot\|_{\mathbf{W^{1,\infty}}}$ norm, as its values are in the space of nonnegative measures. The norm $\|\cdot\|_{\mathbf{BC}^{0,1}}$ in the $\mathbf{BC}^{0,1}$ space is defined in the following way
$$
\|f\|_{\mathbf{BC}^{0,1}}=\sup_{t\in[0,T],\mathbf{v}\in Y}\big(\|f(t,\mathbf{v})\|_X+{\rm \mathbf{Lip_v}}(f(t,\cdot))\big),
$$
where $Y=\mathcal{V}$ or $Y=\mathcal{U}$ and ${\rm \mathbf{Lip_v}}(f(t,\cdot))$ is a Lipschitz constant of $f(t,\cdot)$.
\end{Assumption}
We will write $\|\cdot\|_{\mathbf{BC}}$ instead of $\|\cdot\|_{\mathbf{BC}^{0,1}}$, for the sake of clarity.

\begin{Remark}\label{Lipschitz_semiflow}
It was proven in \cite{Ulikowska2012}, under Assumption \ref{model_functions_measures}, Theorem 2.9, that solutions to (\ref{the_model_measures}) form Lipschitz semiflows $S:[0,T]\times[0,T]\times \mathcal{U}\rightarrow\mathcal{U}$. 
\end{Remark}

\begin{Remark}\label{numerical_solution_semiflow}
Let $\nu_t^m$, $\nu_t^f$, $\nu_t^c$ be defined by (\ref{solution_of_numerical_method}), where $(x_i^m(t), m_i^m(t))$, $(y_j^f(t), m_j^f(t))$, $((x_{ij}^c(t),y_{ij}^c(t)), m_{ij}^c(t))$ is the output of the EBT algorithm (\ref{the_numerical_method_measures}). Then it is easy to show, using Lemma \ref{LemmaFlatMetric}, that the solution of numerical method $\mathbf{v_t}=(\nu_t^m,\nu_t^f,\nu_t^c):[0,T]\rightarrow (\mathcal{U},{\textbf d})$ is a Lipschitz continuous map. The rigorous proof can be seen in \cite[Lemma 4.2]{GOU}.
\end{Remark}

\begin{Remark} \label{initial_accuracy}
It was shown in \cite{JMC}, that for a given Radon measure the initial conditions $u_0^m(x)$ and $u_0^f(y)$ can be approximated in flat metric with an arbitrarily good precision by a combination of Dirac measures defined by (\ref{initial_approximation}). The same reasoning can be applied to the initial condition $u_0^c(\cdot,x,y)$.
\end{Remark}

\begin{Proposition}\label{tangential_ineq}
Let $S:[0,T]\times[0,T]\times E\rightarrow E$ be a Lipschitz semiflow. For every Lipschitz continuous map $[0,T]\ni t\mapsto\nu_t\in E$ the following estimate holds:
\begin{equation*}
d(\nu_t,S(t;0)\nu_0)\leq L\int_{[0,t]}\liminf_{h\rightarrow 0}\frac{d(\nu_{\tau+h},S(h;\tau)\nu_{\tau})}{h}d\tau.
\end{equation*}
\end{Proposition}
The proof of the proposition is similar to the proof of Theorem 2.9 in \cite{bressan_00}.

\begin{Lemma}\label{lem:boundedness_of_density}
Let us consider the equation
\begin{equation}\label{eq:boundedness_of_density}
\begin{array}{rcl}
\partial_t\mu_t+\partial_x\mu_t+\xi(t,\mu)\mu&=&g(t)\delta_{\{0\}},\hspace{1.4cm}(t,x)\in[0,T]\times \bR_+\\
\mu_0^m& \in & \mathcal{M}^+(\bR_+),\\
\end{array}
\end{equation}
where $g(t):=\int_{\bR_+^2}\beta(t,\mu_t)(z)d\mu_t^*(z)$ with the coefficients $\xi$ and $\beta$ satisfying the conditions in Assumptions \ref{model_functions_measures} for the respective coefficients. Here, $\mu_t^*(z)$ is a given curve of measures with finite total mass for all $t\geq 0$ narrowly continuous with respect to time. The weak solution is given by
$$
\mu_t= (\Phi(t;\cdot)\#\mu_0)\exp\left(-\int_0^t\xi(\tau,\cdot -t+\tau){\rm d}\tau\right)+ \tilde \mu_t
$$
with $\Phi(t;x)=x+t$, $\Phi(t;\cdot)\#\mu_0$ refers to the push forward of a measure through the map $\Phi(t;\cdot)$, and
\begin{align*}
\tilde{\mu}_t=\int_0^t\exp\left(-\int_{\tau}^t\xi(s,\cdot+s-t){\rm d}s\right)g(\tau)\delta_{\{t-\tau\}}{\rm d}\tau.
\end{align*}
More precisely, $\tilde \mu_t$ is defined by
$$
\left\langle \varphi,\int_0^t\exp\left(-\int_{\tau}^t\xi(s,\cdot+s-t){\rm d}s\right)g(\tau)\delta_{\{t-\tau\}}{\rm d}\tau\right\rangle=\int_0^t\varphi(t-\tau)\exp\left(-\int_{\tau}^t\xi(s,s-\tau){\rm d}s\right)g(\tau){\rm d}\tau
$$
for all test functions $\varphi$, and its Radon-Nykodym derivative with respect to the Lebesgue measure $\lambda$ on the line is given by the bounded function
\begin{equation*}
f(t,x)=\left\{
\begin{array}{ll}
\displaystyle h(t,t-x)-
\int_x^t\xi(t,t-\tau+x) h(t,\tau-x){\rm d}\tau,& 0\leq x\leq t,\\
0 &x>t.\\
\end{array}
\right. ,
\end{equation*}
that is, $\tilde\mu_t= f(t,\cdot) \lambda$, with $h(t,\tau)=\exp\left(-\int_{\tau}^t\xi(s,s-\tau)
{\rm d}s\right)g(\tau)$.
\end{Lemma}

\begin{proofof}
Let us denote the right-hand side of equation \eqref{eq:boundedness_of_density} as $\omega_t$. It is known that the semigroup of this McKendrick--von Foerster equation can be written in terms of the characteristics of the flow, see \cite{SM2007}. More precisely, defining now 
$\Phi(\tau,t;x)=x+t-\tau$ for any $t,\tau\geq 0$, then the unique solution reads as
\begin{align*}
\mu_t=&\Phi(0,t;\cdot)\#\left(\mu_0\exp\left(-\int_0^t\xi(\tau,\Phi(0,\tau;\cdot)){\rm d}\tau\right)\right)
+\int_0^t\Phi(\tau,t,\cdot)\#\left(\omega_{\tau}\exp\left(-\int_{\tau}^t\xi(s,\Phi(\tau,s;\cdot)){\rm d}s\right)\right){\rm d}\tau
\end{align*}
By substituting $\omega_t=g(t)\delta_{\{0\}}$ and simple algebraic manipulations, we obtain
$$
\mu_t=(\Phi(0,t;\cdot)\#\mu_0)\exp\left(-\int_0^t\xi(\tau,\cdot -t+\tau){\rm d}\tau\right)+\int_0^t\exp\left(-\int_{\tau}^t\xi(s,\cdot -t+s){\rm d}s\right)g(\tau)\delta_{\{t-\tau\}}{\rm d}\tau
$$
as stated above. The final part of the Lemma is obtain by computing directly the Radon-Nikodym derivative using the formula for $\tilde \mu_t$ acting on test functions. Choosing
for sufficiently small $\varepsilon$ the test function defined as 
$$
\varphi_{\varepsilon}=\frac{1}{2\varepsilon}\chi_{(x-\varepsilon,x+\varepsilon)}\,, \quad \mbox{where } \chi_I \mbox{ stands for the characterisctic function of the interval } I\,,
$$
we can check that for $0\leq x\leq t$,
\begin{align*}
\lim_{\varepsilon\to 0}\left\langle \tilde{\mu}_t,\varphi_{\varepsilon}\right\rangle=&\lim_{\varepsilon\to 0}\frac{1}{2\varepsilon}\int_{\max(0,t-x-\varepsilon)}^{\min(0,t-x+\varepsilon)}h(t,\tau){\rm d}\tau 
=\frac{d}{dt}\int_x^tf(t,\tau-x){\rm d}\tau 
=h(t,t-x)+\int_x^t\frac{\partial h}{\partial t}(t,\tau -x){\rm d}\tau \,.
\end{align*}
It is straightforward to check that for $x>t$, the previous computation gives 0, finishing the 
proof of the result.
\end{proofof}

\

The following lemma will be used extensively in the proof of the main theorem on the convergence of the EBT scheme

\begin{Lemma}\label{esti0}
Let  $t\in [\tau, \tau+h)\subset [t_n,t_{n+1})$, and $t\mapsto(\nu_t^m,\nu_t^f,\nu_t^c),\ t\mapsto(\bar{\nu}_t^m,\bar{\nu}_t^f,\bar{\nu}_t^c)\in \mathcal{U}$, where 
$$
\nu_{t}^m:=\sum_{i=B_n}^Jm_i^m(t)\delta_{\{x_i^m(t)\}},\ \nu_{t}^f:=\sum_{j=B_n}^Jm_j^f(t)\delta_{\{y_j^f(t)\}},\ \nu_{t}^c:=\sum_{i,j=B_n}^Jm_{ij}^c(t)\delta_{\{x_{ij}^c(t),y_{ij}^c(t)\}},
$$
and
$$
\bar{\nu}_{t}^m:=\sum_{i=B_n}^J\bar{m}_i^m(t)\delta_{\{x_i^m(t)\}},\ \bar{\nu}_{t}^f:=\sum_{j=B_n}^J\bar{m}_j^f(t)\delta_{\{y_j^f(t)\}},\ \bar{\nu}_{t}^c:=\sum_{i,j=B_n}^J\bar{m}_{ij}^c(t)\delta_{\{x_{ij}^c(t),y_{ij}^c(t)\}},
$$
be certain Lipschitz continuous maps in $d_N$ such, that $(\nu_{\tau}^m,\nu_{\tau}^f,\nu_{\tau}^c)=(\bar{\nu}_{\tau}^m,\bar{\nu}_{\tau}^f,\bar{\nu}_{\tau}^c)$, and 
\begin{equation*}
\begin{split}
C_1=&{\rm Lip}_t(\nu^m)+{\rm Lip}_t(\bar{\nu}^m)+{\rm Lip}_t(\nu^f)+{\rm Lip}_t(\bar{\nu}^f)\\
C_2=&{\rm Lip}_t(\nu^m)+{\rm Lip}_t(\bar{\nu}^m)+{\rm Lip}_t(\nu^f)+{\rm Lip}_t(\bar{\nu}^f)+{\rm Lip}_t(\nu^c)+{\rm Lip}_t(\bar{\nu}^c),
\end{split}
\end{equation*}
where ${\rm Lip}_t(\cdot)$ stands for the Lipschitz constant  with respect to $t$.
Then the following estimates hold for $i,j=B_n,\ldots,J$
\begin{equation}\label{differences_esti}
\begin{split}
\left|\xi^k(t,\nu_t^m,\nu_t^f)(x_i(t))-\xi^k(t,\bar{\nu}_t^m,\bar{\nu}_t^f)(x_i(t))\right|\leq &
C_1h\|\xi^k\|_{\mathbf{BC}},\ k=m,f \\
\left|\beta^k(t,\nu_t^m,\nu_t^f)(x_{ij}(t),y_{ij}(t))-\beta^k(t,\bar{\nu}_t^m,\bar{\nu}_t^f)(x_{ij}(t),y_{ij}(t))\right|\leq &
C_1h\|\beta^k\|_{\mathbf{BC}},\ k=m,f\\
\left|\xi^c(t,\nu_t^m,\nu_t^f,\nu_t^c)(x_{ij}(t),y_{ij}(t))-\xi^c(t,\bar{\nu}_t^m,\bar{\nu}_t^f,\bar{\nu}_t^c)(x_{ij}(t),y_{ij}(t))\right|\leq &
2C_2h\|\xi^c\|_{\mathbf{BC}}.
\end{split}
\end{equation}
Moreover, given functions $m_i(t),\bar{m}_i(t)\in{\rm W^{1,\infty}}([\tau,\tau+h),\bR_+)$ such, that $m_i(\tau)=\bar{m}_i(\tau)$ ($i=B_n,\ldots,J$) and constants
\begin{equation*}
\begin{split}
C_3=&\sup_{t\in[0,T]}\left\{\|\nu_t^m\|_{{\rm TV}}+\|\nu_t^f\|_{{\rm TV}}+\|\nu_t^c\|_{{\rm TV}}+\|\bar\nu_t^m\|_{{\rm TV}}+\|\bar\nu_t^f\|_{{\rm TV}}+\|\bar\nu_t^c\|_{{\rm TV}}\right\},\\
C_4=&{\rm Lip}_t(m_{B_n})+{\rm Lip}_t(\bar{m}_{B_n}),\, C_5=C_1C_3+C_1,\, C_6=2C_2C_3+C_2
\end{split}
\end{equation*}
we can conclude that, for $k=m,f$ and for any $\varphi \in W^{1,\infty}$ with norm less than 1,
\begin{equation}\label{integrals_esti}
\begin{split}
&
\left|\sum_{i=B_n}^J\varphi(x_i(t))[\xi^k(t,\nu_t^m,\nu_t^f)(x_i(t))m_i(t)-\xi^k(t,\bar{\nu}_t^m,\bar{\nu}_t^f)(x_i(t))\bar{m}_i(t)]\right|\leq
C_5 h\|\xi^k\|_{\mathbf{BC}},\\
&
\left|\sum_{i=B_n}^J\varphi(x_i(t))[\beta^k(t,\nu_t^m,\nu_t^f)(x_{ij}(t),y_{ij}(t))m_{ij}(t)-\beta^k(t,\bar{\nu}_t^m,\bar{\nu}_t^f)(x_{ij}(t),y_{ij}(t))\bar{m}_{ij}(t)]\right|\leq
C_5 h\|\beta^k\|_{\mathbf{BC}},\\
&
\Bigg|\sum_{i=B_n}^J\varphi(x_{ij}(t),y_{ij}(t))[\xi^c(t,\nu_t^m,\nu_t^f,\nu_t^c)(x_{ij}(t),y_{ij}(t))m_{ij}(t)-
\xi^c(t,\bar{\nu}_t^m,\bar{\nu}_t^f,\bar{\nu}_t^c)(x_{ij}(t),y_{ij}(t))\bar{m}_{ij}(t)]\Bigg|\leq
C_6 h\|\xi^c\|_{\mathbf{BC}},
\end{split}
\end{equation}
and for $k=m,f$
\begin{equation}\label{integrals_esti_1}
\Bigg|
\int_\tau^{\tau+h}\xi^k(t,\nu_t^m,\nu_t^f)(x_{B_n}(t))m_{B_n}(t)dt-
\int_\tau^{\tau+h}\xi^k(t,\bar{\nu}_t^m,\bar{\nu}_t^f)(x_{B_n}(t))\bar{m}_{B_n}(t)dt\Bigg|\leq (C_4+C_1C_3) h^2\|\xi^k\|_{\mathbf{BC}}.
\end{equation}
\end{Lemma}

\begin{proofof}
We will prove the above estimates for the male function $\xi^m$ only since the same strategy can be used for the remaining functions. To simplify the notation, we avoid the explicit $t$ dependence of $\xi^m$. Let us start with (\ref{differences_esti}). Using the triangular inequality and the assumption 
$(\nu_{\tau}^m,\nu_{\tau}^f,\nu_{\tau}^c)=(\bar{\nu}_{\tau}^m,\bar{\nu}_{\tau}^f,\bar{\nu}_{\tau}^c)$, we have
\begin{align*}
\left|\xi^m(\nu_t^m,\nu_t^f)(x)-\xi^m(\bar{\nu}_t^m,\bar{\nu}_t^f)(x)\right|\leq &
\left|\xi^m(\nu_t^m,\nu_t^f)(x)-\xi^m(\bar{\nu}_t^m,\nu_t^f)(x)\right|+
\left|\xi^m(\bar{\nu}_t^m,\nu_t^f)(x)-\xi^m(\bar{\nu}_t^m,\bar{\nu}_t^f)(x)\right|\\
\leq &\left|\xi^m(\nu_t^m,\nu_t^f)(x)-\xi^m(\nu_\tau^m,\nu_t^f)(x)\right|+
\left|\xi^m(\bar{\nu}_\tau^m,\nu_t^f)(x)-\xi^m(\bar{\nu}_t^m,\nu_t^f)(x)\right|+\\
&\left|\xi^m(\bar{\nu}_t^m,\nu_t^f)(x)-\xi^m(\bar{\nu}_t^m,\nu_\tau^f)(x)\right|+
\left|\xi^m(\bar{\nu}_t^m,\bar{\nu}_\tau^f)(x)-\xi^m(\bar{\nu}_t^m,\bar{\nu}_t^f)(x)\right|
\end{align*}
for all $x\geq 0$ from which
\begin{align*}
\left|\xi^m(\nu_t^m,\nu_t^f)(x)-\xi^m(\bar{\nu}_t^m,\bar{\nu}_t^f)(x)\right|& \leq
\|\xi^m\|_{\mathbf{BC}}\sup_{t\in [\tau, \tau+h)}\left(d_1(\nu_t^m,\nu_\tau^m)+d_1(\bar{\nu}_\tau^m,\bar{\nu}_t^m) + d_1(\nu_t^f,\nu_\tau^f)+
d_1(\bar{\nu}_\tau^f,\bar{\nu}_t^f)\right)
\\
& \leq C_1h\|\xi^m\|_{\mathbf{BC}}
\end{align*}
for all $x\geq 0$.

Let us move to the proof of (\ref{integrals_esti}) now, for any $\varphi \in W^{1,\infty}$ with norm less than 1 and for all $x\geq0$, we have
\begin{align*}
\left|\sum_{i={B_n}}^J\varphi(x)[\xi^m(\nu_t^m,\nu_t^f)(x)m_i(t)-\xi^m(\bar{\nu}_t^m,\bar{\nu}_t^f)(x)\bar{m}_i(t)]\right|\leq &
\left|\sum_{i=B_n}^J\varphi(x)[\xi^m(\nu_t^m,\nu_t^f)(x)m_i(t)-\xi^m(\bar{\nu}_t^m,\bar{\nu}_t^f)(x)m_i(t)]\right|+\\
&
\left|\sum_{i=B_n}^J\varphi(x)\xi^m(\bar{\nu}_t^m,\bar{\nu}_t^f)(x)[m_i(t)-\bar{m}_i(t)]\right|
\end{align*}
Using inequalities (\ref{differences_esti}) and the fact that $m_i(t)\geq0$ for $i=B_n,\ldots,J$, the first term of the right hand side of the above inequality can be estimated as
$$
\left|\sum_{i=B_n}^J\varphi(x)[\xi^m(\nu_t^m,\nu_t^f)(x)m_i(t)-\xi^m(\bar{\nu}_t^m,\bar{\nu}_t^f)(x)m_i(t)]\right|\leq
C_1h\|\xi^m\|_{\mathbf{BC}}\sum_{i=B_n}^J\left|\varphi(x)\right| m_i(t)\leq
C_1C_3h\|\xi^m\|_{\mathbf{BC}}\,.
$$
To estimate the second term of the right hand side, we proceed as
$$
\left|\sum_{i=B_n}^J\varphi(x)\xi^m(\bar{\nu}_t^m,\bar{\nu}_t^f)(x)[m_i(t)-\bar{m}_i(t)]\right|\leq\|\xi^m(\nu_t^m,\nu_t^f)\|_{L^\infty} d_1(\nu_t^m,\bar{\nu}_t^m)\leq C_1 h\|\xi^m\|_{\mathbf{BC}},
$$
for any $\varphi \in W^{1,\infty}$ with norm less than 1 and for all $x\geq0$.
Here, we used that 
$$
d_1(\nu_t^m,\bar{\nu}_t^m)=
\sum_{i=B_n}^J|m_i(t)-\bar{m}_i(t)|\leq 
d_1(\nu_t^m,\nu_\tau^m) + 
d_1(\bar{\nu}_\tau^m,\bar{\nu}_t^m)
$$
since both measures are combinations of Dirac Deltas at the same locations and assumption $(\nu_{\tau}^m,\nu_{\tau}^f,\nu_{\tau}^c)=(\bar{\nu}_{\tau}^m,\bar{\nu}_{\tau}^f,\bar{\nu}_{\tau}^c)$.
The proof for (\ref{integrals_esti_1}) is easier as
\begin{align*}
&\left|\int_\tau^{\tau+h}\left(\xi^m(t,\nu_t^m,\nu_t^f)(x_{B_n}(t))m_{B_n}^m(t)-\xi^m(t,\bar{\nu}_t^m,\bar{\nu}_t^f)(x_{B_n}(t))\bar{m}_{B_n}^m(t)\right)dt\right|\leq \\
&\quad\int_\tau^{\tau+h}\!\!\left|\xi^m(t,\nu_t^m,\nu_t^f)(x_{B_n}(t))-\xi^m(t,\bar{\nu}_t^m,\bar{\nu}_t^f)(x_{B_n}(t))\right|\bar{m}_{B_n}^m(t)dt+\int_\tau^{\tau+h}\!\!\!\!\!\xi^m(t,\nu_t^m,\nu_t^f)(x_{B_n}(t))\left|m_{B_n}^m(t)-\bar{m}_{B_n}^m(t)\right|dt
\\
& \quad \leq C_1 h\|\xi^m\|_{\mathbf{BC}}\int_\tau^{\tau+h}\bar{m}_{B_n}^m(t)dt+\|\xi^m\|_{\mathbf{BC}}\int_\tau^{\tau+h}h({\rm Lip}_t(m_{B_n})+{\rm Lip}_t(\bar{m}_{B_n}))dt\leq \left(C_1C_3+C_4\right) h^2\|\xi^m\|_{\mathbf{BC}}.
\end{align*}
\end{proofof}


\section{Convergence}
\setcounter{equation}{0}

\begin{Theorem}\label{convergence}
Let $\mathbf{u}_t=(\mu_t^m,\mu_t^f,\mu_t^c):[0,T]\rightarrow \mathcal{U}$ be solution of (\ref{the_model_measures}), and $\mathbf{v_0}=(\nu_0^m,\nu_0^f,\nu _0^c) \in \mathcal{U}$ be an approximation of $\mathbf{u_0}=(\mu_0^m,\mu_0^f,\mu _0^c)$ given by formulas (\ref{initial_approximation}), (\ref{the_numerical_method_initial_2}) and (\ref{the_numerical_method_initial_3}) with an error
$$
\varepsilon_0=\mathbf{d}(\mathbf{u_0},\mathbf{v_0}).
$$
Let $\mathbf{v_t}=(\nu_t^m,\nu_t^f,\nu _t^c):[0,T]\rightarrow \mathcal{U}$ be the output of the numerical method defined in Section \ref{EBT_measures} with the initial condition $\mathbf{v_0}$. Then, there exists such a constant $C$, that
\begin{equation}\label{global}
\mathbf{d}(\mathbf{\mathbf{u}_t},\mathbf{\mathbf{v}_t})\leq\varepsilon_0+C T\Delta t,\ t\leq T.
\end{equation}
where $\Delta t=\max_{n=0,\ldots {N_T}-1} |t_{n+1}-t_n|$,  $t_n<T$ for $n=0,\ldots,{N_T}$ and $\Delta t\leq a_0$, where $a_0$ satisfies \eqref{newbornsdonotmarry}.
\end{Theorem}

\begin{proofof}
As it was mentioned in Remark \ref{Lipschitz_semiflow} and Remark \ref{numerical_solution_semiflow}, the solution of the numerical method defined in Section \ref{EBT_measures} ${\textbf v}_t$ is a Lipschitz continues map, and the problem (\ref{the_model_measures}) generates a Lipschitz semiflow $S$ in metric space $(\mathcal{U},{\mathbf d})$. Thus, to estimate the difference between ${\textbf v}_t$ and ${\textbf u}_t=S(t;0){\textbf u}_0$, we are entitled to use Proposition \ref{tangential_ineq} and to consider the problem locally in time.

Without loss of generality we can assume, that interval $(\tau,\tau+h]$ does not contain internalisation point time, that is, there exists $n$ such that $(\tau,\tau+h]\subset [t_n,t_{n+1})$. A remark on this assumption will be made just after the proof.

According to Proposition \ref{tangential_ineq}, we wish to estimate the following distance
\begin{equation}\label{whole_distance}
{\mathbf d}({\textbf v}_{\tau+h},S(h;\tau){\textbf v}_\tau)=d_1(\nu^m_{\tau+h},S^m(h;\tau)\nu^m_\tau)+d_1(\nu^f_{\tau+h},S^f(h;\tau)\nu^f_\tau)
+d_2(\nu^c_{\tau+h},S^c(h;\tau)\nu^c_\tau)
\end{equation}
showing that it is bounded by $C(t_{n+1}-t_n) \Delta t$ with $C$ depending only on $T$ and the hypotheses on the functions involved in (\ref{the_model_measures}) given by Assumption \ref{model_functions_measures}.
It is straightforward from this local estimate using the properties of the $d_N$-semiflow in Definition \ref{semiflow} and the triangular inequality to deduce the global estimate \eqref{global} by induction.

Let us recall, that components of the EBT ODE output ${\textbf v}_{t}$, $t\in[t_n,t_{n+1})$, are defined as linear combinations of $J-B_n+1$ Dirac measures in case of males/females, and $(J-B_n+1)^2$ in case of couples
$$
\nu_t^m:=\sum_{i=B_n}^Jm_i^m(t)\delta_{\{x_i^m(t)\}},\ \nu_t^f:=\sum_{j=B_n}^Jm_j^f(t)\delta_{\{y_j^f(t)\}},\ 
\nu_{t}^c:=\sum_{i,j=B_n}^Jm_{ij}^c(t)\delta_{\{x_{ij}^c(t),y_{ij}^c(t)\}}
$$  
Let $\bar{{\mathbf v}}_t=(\bar{\nu}_{t}^m,\bar{\nu}_{t}^f,\bar{\nu}_{t}^c)$ be a solution to (\ref{the_model_measures}) on the interval $[\tau,\tau+h]$ with initial condition ${\textbf v}_\tau$
$$
\bar{\nu}^m_{t}:=S^m(t-\tau;\tau)\nu^m_{\tau},\hspace{1cm} \bar{\nu}^f_{t}:=S^f(t-\tau;\tau)\nu^f_{\tau},\hspace{0.5cm} {\rm and}\ \hspace{0.5cm} \bar{\nu}^c_{t}:=S^c(t-\tau;\tau)\nu^c_{\tau},
$$

The semiflow $S=(S^m,S^f,S^c)$ generates the following solution
$$
\bar{\nu}^m_t=f^m(t,\cdot)\lambda^1+\sum_{i=B_n}^J \bar{m}_i^m(t)\delta_{\{\bar{x}_i^m(t)\}},\ \ \
\bar{\nu}^f_{t}=f^f(t,\cdot)\lambda^1+\sum_{j=B_n}^J \bar{m}_j^f(t)\delta_{\{\bar{y}_j^f(t)\}},\ \ \
\bar{\nu}^c_{t}=\sum_{i,j=B_n}^J \bar{m}_{ij}^c(t)\delta_{\{\bar{x}_{ij}^c(t),\bar{y}_{ij}^c(t)\}},
$$
where $\lambda^1$ is a Lebesgue measure, and $f^m$ and $f^f$ are densities arising in the boundary cohorts, such that 
${\rm supp} f^k(t,\cdot)\subset(0,t-\tau),\ k=f,m$. Observe, that $\bar{\nu}^c_t$ does not have an absolutely continuous measure $f^c(t,)\lambda^2$ because newborns do not form couples. This is guaranteed by marriage function \eqref{spec_marriage_function} together with assumptions \eqref{newbornsdonotmarry} and $\Delta t\leq a_0$. 

Let us write $s=\tau+h$ and notice that the first component of (\ref{whole_distance}) can be initially estimated in the following way
\begin{equation}\label{d1}
\begin{split}
d_1(\bar{\nu}_{s}^m,\nu_{s}^m)
=&\,
d_1\left(f^m(s,\cdot)\lambda^1+\sum_{i=B_n}^J\bar{m}^m_{i}(s)\delta_{\{\bar{x}^m_i(s)\}},\sum_{i=B_n}^Jm^m_{i}(s)\delta_{\{x^m_i(s)\}}\right)\\
\leq&\,
d_1\left(f^m(s,\cdot)\lambda^1,\bar{p}_{B_n}^m(s)\delta_{\{\bar{x}^m_{B_n}(s)\}}\right)+d_1\left(\bar{p}^m_{B_n}(s)\delta_{\{\bar{x}^m_{B_n}(s)\}}+\sum_{i=B_n}^J\bar{m}^m_{i}\delta_{\{\bar{x}^m_i(s)\}},\sum_{i=B_n}^Jm^m_{i}
\delta_{\{x^m_i(s)\}}\right)\\
\leq&\,
d_1\left(f^m(s,\cdot)\lambda^1,\bar{p}^m_{B_n}(s)\delta_{\{\bar{x}^m_{B_n}(s)\}}\right)+
d_1\left(\left(\bar{p}^m_{B_n}(s)+\bar{m}^m_{B_n}(s)\right)\delta_{\{\bar{x}^m_{B_n}(s)\}},m^m_{B_n}(s)\delta_{\{x^m_{B_n}(s)\}}\right)\\
&+
\sum_{i=B_n+1}^Jd_1\left(\bar{m}^m_{i}(s)\delta_{\{\bar{x}^m_i(s)\}},m^m_{i}(s)\delta_{\{x^m_i(s)\}}\right)\\
\end{split}\end{equation}
where $\bar{p}_{B_n}^m(s)$ is a mass generated by the absolutely continuous measures $f^m(t,\cdot)\lambda^1$, i.e. $\bar{p}_{B_n}^m(s)=\int_{\bR_+}f^m(s,x)dx$. The previous estimate uses that the flat metric is in fact a norm. The second component of (\ref{whole_distance}) concerning the female population is treated analogously, while the third component will follow easily from the expression
\begin{equation}\label{d2}
d_2(\bar{\nu}_{s}^c,\nu_{s}^c)=d_2\left(\sum_{i,j=B_n}^J \bar{m}_{ij}^c(s)\delta_{\{\bar{x}_{ij}^c(s),\bar{y}_{ij}^c(s)\}},\sum_{i,j=B_n}^J m_{ij}^c(s)\delta_{\{x_{ij}^c(s),y_{ij}^c(s)\}}\right).
\end{equation}

To handle the estimates (\ref{d1}) and (\ref{d2}) it is necessary to find the locations $\bar{x}_i^m(s)$, $\bar{y}_j^f(s)$, $(\bar{x}_{ij}^c(s),\bar{y}_{ij}^c(s))$ and masses $\bar{m}_i^m(s)$, $\bar{m}_j^f(s)$, $\bar{m}_{ij}^c(s)$, $\bar{p}_{B_n}^m(s)$, $\bar{p}_{B_n}^f(s)$ for $i,j=B_n,\ldots,J$, which are generated by the semiflow $S$.

\begin{enumerate}
\item[\textbf{I}]
\textbf{- Locations}. Let us start with male population and observe that locations $x_i^m(s)$ and $\bar{x}_i^m(s)$, $i=B_n,\ldots,J$ are equal in $[\tau,\tau+h]$, as they are governed by the same rule ($t\in[\tau,\tau+h]$)
$$
\begin{array}{rcl}
\frac{d}{dt}x_i^m(t)&=&1,\\
\frac{d}{dt}\bar{x}_i^m(t)&=&1,\\
\bar{x}_i^m(t)&=&x_i^m(t),\ i=B_n,\ldots,J.
\end{array}
$$
Using similar argument for female and couple populations we end up with the following equalities 
\begin{equation*}
\begin{split}
x_i^m(t)&=\bar{x}_i^m(t),\\
y_i^f(t)&=\bar{y}_i^f(t),\\
(x_{ij}^c(t),y_{ij}^c(t))&=(\bar{x}_{ij}^c(t),\bar{y}_{ij}^c(t)),\ {\rm for}\ i,j=B_n,\ldots,J.
\end{split}
\end{equation*}
\item[\textbf{II}]
\textbf{- Masses}. The following formulas for masses at time $s=\tau+h$,
\begin{equation*}
\begin{split}
\bar{m}_i^m(s)=\ &m_i^m(\tau)-\int_\tau^{\tau+h}\xi^m(t,\bar{\nu}_t^m,\bar{\nu}_t^f)(x_i^m(t))\bar{m}_i^m(t)dt,\\
\bar{p}^m_{B_n}(s)=&\int_{\tau}^{\tau+h}\int_{\bR_+^2}\beta^m(t,\bar{\nu}_t^m,\bar{\nu}_t^f)(z)d\bar{\nu}_t^c(z)dt-\int_{\tau}^{\tau+h}\int_{\bR_+}\xi^m(t,\bar{\nu}_t^m,\bar{\nu}_t^f)f^m(t,x)dxdt,\\
\bar{m}_j^f(s)=\ &m_j^f(\tau)-\int_\tau^{\tau+h}\xi^f(t,\bar{\nu}_t^m,\bar{\nu}_t^f)(y_j^f(t))\bar{m}_j^f(t)dt,\\
\bar{p}^f_{B_n}(s)=&\int_{\tau}^{\tau+h}\int_{\bR_+^2}\beta^f(t,\bar{\nu}_t^m,\bar{\nu}_t^f)(z)d\bar{\nu}_t^c(z)dt-\int_{\tau}^{\tau+h}\int_{\bR_+}\xi^f(t,\bar{\nu}_t^m,\bar{\nu}_t^f)f^f(t,y)dydt,\\
\bar{m}^c_{ij}(s)=\ &m_{ij}^c(\tau)-\int_\tau^{\tau+h}\xi^m(t,\bar{\nu}_t^m,\bar{\nu}_t^f,\bar{\nu}_t^c)(x_{ij}^c(t),y_{ij}^c(t))\bar{m}_{ij}^c(t)dt
+\int_\tau^{\tau+h}\frac{\tilde{N}_{ij}(t)}{\tilde{D}_{ij}(t)}dt,
\end{split}
\end{equation*}
will be derived using three different types of test functions: $\varphi_\varepsilon(t,x)$, $\varphi_\varepsilon^i(t,x)$ and $\varphi_\varepsilon^{ij}(t,x,y),\ i,j=B_n,\ldots,J$. In all those test functions parameter should be chosen in such a way,  that the support of each test function does intersect with the domain of only  one, particular location function $x_{i}^m(t)$, $(x_{ij}^c,y_{ij}^c)(t)$, $t\in[\tau,\tau+h)$ and $i,j=B_n,\ldots,J$. Such a choice is possible due to the regularity of cohorts boundaries, which are straight parallel lines.

\begin{enumerate}
\item[\textbf{(a)}]
$\bar{m}_i^m(s)$, $\bar{m}_j^f(s)$, $i,j=B_n,\ldots,J$.\\
To derive the evolution of the male population mass, which is generated by the semiflow in the $i$-th internal cohorts, we use the following test functions $\varphi^i_{\varepsilon}\in(\mathbf{C^1}\cap\mathbf{W^{1,\infty}})([\tau,\tau+h]\times\bR_+;\bR)$
$$
\varphi^i_{\varepsilon}(t,x)=\left\{
\begin{array}{l}
1\ {\rm if}\ x\in[x_{i}^m(\tau)-\varepsilon,x_{i}^m(\tau+h)+\varepsilon],\\[2mm]
\frac{x-(x_{i}^m(\tau)-2\varepsilon)}{\varepsilon}\ {\rm if}\ x\in[x_{i}^m(\tau)-2\varepsilon,x_{i}^m(\tau)-\varepsilon],\\[2mm]
\frac{-x+(x_{i}^m(\tau+h)+2\varepsilon)}{\varepsilon}\ {\rm if}\ x\in[x_{i}^m(\tau+h)+\varepsilon,x_{i}^m(\tau+h)+2\varepsilon],\\[2mm]
0,\ {\rm otherwise},
\end{array}
\right.
$$
According to Definition \ref{weak_solution}, if the measure $(\bar{\nu}_t^m,\bar{\nu}_t^f,\bar{\nu} _t^c)$ is a weak solution to (\ref{the_model_measures}) on time interval $[\tau,\tau+h]$ (not $[0,T]$, as we investigate the equation locally in time) then the following equality holds
$$
\int_{\tau}^{\tau+h}\int_{\bR_+}\left(\partial_t\varphi^i_{\varepsilon}(t,x)+\partial_x\varphi^i_{\varepsilon}(t,x)-\xi^m(t,\bar{\nu}_t^m,\bar{\nu}_t^f)\varphi^i_{\varepsilon}(t,x)\right)d\bar{\nu}_t^m(x)dt
$$
$$
+\int_{\tau}^{\tau+h}\varphi^i_{\varepsilon}(t,0)\int_{\bR_+^2}\beta^m(t,\bar{\nu}_t^m,\bar{\nu}_t^f)(z)d\bar{\nu}_t^c(z)dt
$$
$$
=\int_{\bR_+}\varphi^i_{\varepsilon}(\tau+h,x)d\bar{\nu}_{\tau+h}^m(x)-\int_{\bR_+}\varphi^i_{\varepsilon}(\tau,x)d\bar{\nu}_{\tau}^m(x).
$$

The following integrals vanish: $\int_{\tau}^{\tau+h}\int_{\bR_+}\partial_t\varphi^i_{\varepsilon}(t,x)d\bar{\nu}_t^m(x)dt$, because the test functions do not depend on $t$, $\int_{\tau}^{\tau+h}\int_{\bR_+}\partial_x\varphi^i_{\varepsilon}(t,x)d\bar{\nu}_t^m(x)dt$, because of measure $\bar{\nu}_t^m(x)$ and $\int_{\tau}^{\tau+h}\varphi^i_{\varepsilon}(t,0)\int_{\bR_+^2}\beta^m(t,\bar{\nu}_t^m,\bar{\nu}_t^f)(z)d\bar{\nu}_t^c(z)dt$ because of the support of test functions. Passing with $\varepsilon$ to $0$, and using dominated  convergence theorem we finally obtain the formula for the desired coefficient $\bar{m}_i(s)$, $i=B_n,\ldots,J$:
$$
\bar{m}_i^m(s)=m_i^m(\tau)-\int_\tau^{\tau+h}\xi^m(t,\bar{\nu}_t^m,\bar{\nu}_t^f)(x_i^m(t))\bar{m}_i^m(t)dt.
$$
Formulas for female masses $m_j^f(s)$, $j=1,\ldots,J$ are derived using the same type of test function.

\item[\textbf{(b)}]
${\bar{p}_{B_n}^m(s)}$, ${\bar{p}_{B_n}^f(s)}$. \\
In order to find $\bar{p}_{B_n}^m(s)=\int_{\bR_+}f^m(s,x)dx$ 
we define the test function $\varphi_{\varepsilon}\in(\mathbf{C^1}\cap\mathbf{W^{1,\infty}})([\tau,\tau+h]\times\bR_+;\bR)$ in the following way
$$
\varphi_{\varepsilon}(t,x)=\left\{
\begin{array}{l}
1\ {\rm if}\ x\in[0,h+\varepsilon],\\[2mm]
\frac{-x+(h+2\varepsilon)}{\varepsilon}\ {\rm if}\ x\in[h+\varepsilon,h+2\varepsilon],\\[2mm]
0,\ {\rm otherwise},
\end{array}
\right.
$$
Like previously, we apply the above test function to the Definition \ref{weak_solution} and observe that $\partial_t\varphi_{\varepsilon}=0$, that the ${\rm supp}\ \partial_x\varphi_{\varepsilon}(t,\cdot)\cap{\rm supp}\ f^m(t,\cdot)=\emptyset$, so the first two integrals vanish. Passing with $\varepsilon$ to $0$, we obtain
$$
-\int_{\tau}^{\tau+h}\int_{\bR_+}\xi^m(t,\bar{\nu}_t^m,\bar{\nu}_t^f)f^m(t,x)dxdt+\int_{\tau}^{\tau+h}\int_{\bR_+^2}\beta^m(t,\bar{\nu}_t^m,\bar{\nu}_t^f)(z)d\bar{\nu}_t^c(z)dt
=\int_{\bR_+}f^m(\tau+h,x)dx,
$$
and conclude that
$$
\bar{p}^m_{B_n}(\tau+h)=-\int_{\tau}^{\tau+h}\int_{\bR_+}\xi^m(t,\bar{\nu}_t^m,\bar{\nu}_t^f)f^m(t,x)dxdt+\int_{\tau}^{\tau+h}\int_{\bR_+^2}\beta^m(t,\bar{\nu}_t^m,\bar{\nu}_t^f)(z)d\bar{\nu}_t^c(z)dt.
$$
The same reasoning should be applied to female case.

\item[\textbf{(c)}]
${\bar{m}_{ij}^c(s)}$, $i,j=B_n,\ldots,J$.\\
We will find the evolution of the masses for couples $\bar{m}_{ij}^c(c)$, $i,j=B_n,\ldots,J$, using the following test function $\varphi^{ij}_{\varepsilon}\in(\mathbf{C^1}\cap\mathbf{W^{1,\infty}})([\tau,\tau+h]\times\bR_+\times\bR_+;\bR)$:
$$
\varphi^{ij}_{\varepsilon}(t,x,y)=\left\{
\begin{array}{lll}
1&{\rm if}\ x\in[a,b],&y\in[c,d],\\[2mm]
\frac{x-(a-\varepsilon)}{\varepsilon}&{\rm if}\ x\in[a-\varepsilon,a],&y\in[c,d],\\[2mm]
\frac{-x+(b+\varepsilon)}{\varepsilon}& {\rm if}\ x\in[b,b+\varepsilon],&y\in[c,d],\\[2mm]
\frac{y-(c-\varepsilon)}{\varepsilon}& {\rm if}\ x\in[a,b],&y\in[c-\varepsilon,c],\\[2mm]
\frac{-y+(d+\varepsilon)}{\varepsilon}& {\rm if}\ x\in[a,b],&y\in[d,d+\varepsilon],\\[2mm]
\frac{x+y-(a+c-\varepsilon)}{\varepsilon}& {\rm if}\ x\in[a-\varepsilon,a],&y\in[-x+a+c-\varepsilon,c],\\[2mm]
\frac{x-y-a+\varepsilon+d}{\varepsilon}& {\rm if}\ x\in[a-\varepsilon,a],&y\in[d,x-a+\varepsilon+d],\\[2mm]
\frac{-x+y+b-c+\varepsilon}{\varepsilon}& {\rm if}\ x\in[b,b+\varepsilon],&y\in[x-b+c-\varepsilon,c],\\[2mm]
\frac{-x-y+b+d+\varepsilon}{\varepsilon}&{\rm if}\ x\in[b,b+\varepsilon],&y\in[d,-x+b+d+\varepsilon],\\[2mm]
0,& {\rm otherwise},&
\end{array}
\right.
$$
where
\begin{equation*}
\begin{split}
a=&x_{ij}^c(\tau)-\varepsilon, \hspace{2cm} b=x_{ij}^c(\tau+h)+\varepsilon,\\
c=&y_{ij}^c(\tau)-\varepsilon, \hspace{2cm} d=y_{ij}^c(\tau+h)+\varepsilon,
\end{split}
\end{equation*}

Again, according to Definition \ref{weak_solution}, if the measure $(\bar{\nu}_t^m,\bar{\nu}_t^f,\bar{\nu} _t^c)$ is a weak solution to (\ref{the_model_measures}) on time interval $[\tau,\tau+h]$ then the following equality holds
\begin{align*}
\int_\tau^{\tau+h}\int_{\bR_+}&\left(\partial_t\varphi^{ij}_{\varepsilon}(t,z)+\partial_x\varphi^{ij}_{\varepsilon}(t,z)+\partial_y\varphi^{ij}_{\varepsilon}(t,z)-\xi^c(t,\bar{\nu}_t^m,\bar{\nu}_t^f,\bar{\nu}_t^c)\varphi^{ij}_{\varepsilon}(t,z)\right)d\bar{\nu}_t^c(z)dt\\
&+\int_\tau^{\tau+h}\int_{\bR_+^2}\varphi^{ij}_{\varepsilon}(t,z)d\mathcal{T}(t,\bar{\nu}_t^m,\bar{\nu}_t^f,\bar{\nu}_t^c)(z)dt
=\int_{\bR_+}\varphi^{ij}_{\varepsilon}({\tau+h},z)d\bar{\nu}_{\tau+h}^c(z)-\int_{\bR_+}\varphi^{ij}_{\varepsilon}(\tau,z)d\bar{\nu}_\tau^c(z),
\end{align*}
which, after following similar arguments to the male population leads to
\begin{align*}
\bar{m}^c_{ij}(s)=&\,m_{ij}^c(\tau)-\int_\tau^{\tau+h}\xi^m(t,\bar{\nu}_t^m,\bar{\nu}_t^f,\bar{\nu}_t^c)(x_{ij}^c(t),y_{ij}^c(t))\bar{m}_{ij}^c(t)dt\\
&+\lim_{\varepsilon\rightarrow0}\int_\tau^{\tau+h}\int_{\bR_+^2}\varphi^{ij}_{\varepsilon}(t,z)d\mathcal{T}(t,\bar{\nu}_t^m,\bar{\nu}_t^f,\bar{\nu}_t^c)(z)dt.
\end{align*}

Let us observe now, that the measure $\mathcal{T}$, analogously to Definitions (\ref{newlyweds}), (\ref{singles}) and (\ref{spec_marriage_function_measures}), is given by 
$$
\mathcal{T}(t,\bar{\nu}_t^m,\bar{\nu}_t^f,\bar{\nu}_t^c))=\frac{\Theta(x,y)h(x)g(y)}{\gamma+\int_0^\infty h(z)d\bar{s}_t^m(z)+\int_0^\infty g(w)d\bar{s}_t^f(w)}(\bar{s}_t^m\otimes \bar{s}_t^f),
$$
where
$$
\bar{s}_t^m(B)=(\bar{\mu}_t^m-\bar{\sigma}_t^m)(B\times\bR_+),\  \bar{s}_t^f(B)=(\bar{\mu}_t^f-\bar{\sigma}_t^f)(B\times\bR_+)
$$
and
$$
\bar{\sigma}_t^m(B)=(\bar{\mu}_t^c)(B\times\bR_+),\  \bar{\sigma}_t^f(B)=\bar{\mu}_t^c(\bR_+\times B).
$$
The fact that for every fixed $i,j$ ${\rm supp}\ \varphi^{ij}_{\varepsilon}(t,\cdot)\cap{\rm supp}\ f^m(t,\cdot)=\emptyset$ and ${\rm supp}\ \varphi^{ij}_{\varepsilon}(t,\cdot)\cap{\rm supp}\ f^f(t,\cdot)=\emptyset$ yields
$$
\lim_{\varepsilon\rightarrow0}\int_\tau^{\tau+h}\int_{\bR_+^2}\varphi^{ij}_{\varepsilon}(t,z)d\mathcal{T}(t,\bar{\nu}_t^m,\bar{\nu}_t^f,\bar{\nu}_t^c)(z)dt=\int_\tau^{\tau+h}\frac{\tilde{N}_{ij}(t)}{\tilde{D}_{ij}(t)}dt,
$$
where
$$
\frac{\tilde{N}_{ij}(t)}{\tilde{D}_{ij}(t)}=\frac{\Theta(x_{ij}^c(t),y_{ij}^c(t))h(x_{ij}^c(t))g(y_{ij}^c(t))\left(\bar{m}_{i}^m(t)-\sum_{w=B_n}^J\bar{m}_{iw}^c(t)\right)\left(\bar{m}_{j}^f(t)-\sum_{v=B_n}^J\bar{m}_{vj}^c(t)\right)}
{\gamma+\sum_{v=B_n}^Jh(x_{vj}^c(t))\left(\bar{m}_{v}^m(t)-\sum_{w=B_n}^J\bar{m}_{vw}^c(t)\right)+
\sum_{w=B_n}^J g(y_{iw}^c(t))\left(\bar{m}_{w}^f(t)-\sum_{v=B_n}^J\bar{m}_{vw}^c(t)\right)}.
$$
\end{enumerate}
\end{enumerate}

Having desired formulas for locations and masses derived, we can proceed with inequalities (\ref{d1}) and (\ref{d2}). Constants $C_1,C_2$ and $C_3$ are defined exactly like in Lemma \ref{esti0}, while $C_4$ is defined similarly -- but with respect to the underlying mass $m$. Using Lemmas \ref{LemmaFlatMetric} and \ref{esti0}, we will show the following estimates for each term:
\begin{enumerate}
\item ${d_1\left(f^m(s,\cdot)\lambda^1,\bar{p}^m_{B_n}(s)\delta_{\{\bar{x}^m_{B_n}(s)\}}\right)=\mathcal{O}(h)\Delta t}$\\
Let us observe, that
\begin{align*}
d_1&\left(f^m(s,\cdot)\lambda^1,\bar{p}^m_{B_n}(s)\delta_{\{\bar{x}^m_{B_n}(s)\}}\right)\leq \bar{p}^m_{B_n}(s)\bar{x}^m_{B_n}(s)
\\
&\leq\left|-\int_{\tau}^{\tau+h}\int_0^{t-\tau}\xi^m(t,\bar{\nu}_t^m,\bar{\nu}_t^f)f^m(t,x)dxdt+\int_{\tau}^{\tau+h}\int_{\bR_+^2}\beta^m(t,\bar{\nu}_t^m,\bar{\nu}_t^f)(z)d\bar{\nu}_t^c(z)dt\right| \bar{x}_{B_n}^m(s)
\\
&\leq (C_7 h^2\|\xi^m\|_{\mathbf{BC}}+C_8 h\|\beta^m\|_{\mathbf{BC}})\Delta t,
\end{align*}
because $\bar{x}^m_{B_n}(t)=t-\tau\leq\Delta t$ and $f^m$ is bounded due to Lemma \ref{lem:boundedness_of_density} since $\xi^m$ and $\beta^m$ are absolutely continuous and $\bar{\nu}^c$ is finite on $[0,T]$. Here, we defined 
$C_7=\max_{t\in[\tau,\tau+h],\ x\in[0,t-\tau]}f^m(t,x)$
and $C_8=\int_{\bR_+^2}d\bar{\nu}_t^c(z)$.

\item ${d_1\left(\left(\bar{p}^m_{B_n}(s)+\bar{m}^m_{B_n}(s)\right)\delta_{\{\bar{x}^m_{B_n}(s)\}},m^m_{B_n}(s)\delta_{\{x^m_{B_n}(s)\}}\right)=\mathcal{O}(h^2)}$

\begin{align*}
d_1&\left(\left(\bar{p}^m_{B_n}(s)+\bar{m}^m_{B_n}(s)\right)\delta_{\{\bar{x}^m_{B_n}(s)\}},m^m_{B_n}(s)\delta_{\{x^m_{B_n}(s)\}}\right)\leq 
|\bar{p}^m_{B_n}(s)+\bar{m}^m_{B_n}(s)-m^m_{B_n}(s)|=
\\
&\leq
\Bigg|\int_{\tau}^{\tau+h}\int_0^{t-\tau}\xi^m(t,\bar{\nu}_t^m,\bar{\nu}_t^f)f^m(t,x)dxdt\Bigg|
\\
& \quad + \Bigg|
\int_\tau^{\tau+h}\xi^m(t,\nu_t^m,\nu_t^f)(x_{B_n}^m(t))m_{B_n}^m(t)dt-
\int_\tau^{\tau+h}\xi^m(t,\bar{\nu}_t^m,\bar{\nu}_t^f)(x_{B_n}^m(t))\bar{m}_{B_n}^m(t)dt\Bigg|
\\
&\quad +\Bigg|
\int_\tau^{\tau+h}\!\!\!\!\sum_{i,j={B_n}}^J\beta^m(t,\bar{\nu}_t^m,\bar{\nu}_t^f)(x_{ij}^c(t),y_{ij}^c(t))\bar{m}_{ij}^c(t)dt-\!
\int_\tau^{\tau+h}\!\!\!\!\sum_{i,j={B_n}}^J\beta^m(t,\nu_t^m,\nu_t^f)(x_{ij}^c(t),y_{ij}^c(t))m_{ij}^c(t)dt
\Bigg|
\\
&\leq C_5 h^2 \|\xi^m\|_{\mathbf{BC}}+\left(C_1C_3+C_4\right) h^2\|\xi^k\|_{\mathbf{BC}}+\left(C_1C_3+C_1\right) h^2\|\beta^k\|_{\mathbf{BC}},
\end{align*}
where the definitions of the masses $\bar{p}^m_{B_n}$, $\bar{m}^m_{B_n}$, and $m^m_{B_n}$, and $\bar{\nu}_t^c$ were used together with similar estimates as above in the first item.

\item ${d_1\left(\sum_{i=B_n+1}^J\bar{m}^m_{i}(s)\delta_{\{\bar{x}^m_i(s)\}},\sum_{i=B_n+1}^Jm^m_{i}(s)\delta_{\{x^m_i(s)\}}\right)=\mathcal{O}(h^2)}$

Using that the characteristics verify $x^m_i(s)=\bar{x}^m_i(s)$ and \eqref{integrals_esti}, we get
\begin{align*}
d_1&\left(\sum_{i=B_n+1}^J\bar{m}^m_{i}(s)\delta_{\{\bar{x}^m_i(s)\}},\sum_{i=B_n+1}^Jm^m_{i}(s)\delta_{\{x^m_i(s)\}}\right)
\\
&\leq \sup_{\|\varphi\|_{W^{1,\infty}}\leq1}\left|\int_\tau^{\tau+h}\sum_{i=1}^N\varphi(x_i(t))[\xi^m(t,\nu_t^m,\nu_t^f)(x_i(t))m_i(t)-\xi^m(t,\bar{\nu}_t^m,\bar{\nu}_t^f)(x_i(t))\bar{m}_i(t)]dt\right|\\
&\leq
C_5 h^2\|\xi^m\|_{\mathbf{BC}}.
\end{align*}

\item ${d_2\left(\sum_{i,j=B_n}^J \bar{m}_{ij}^c(s)\delta_{\{\bar{x}_{ij}^c(s),\bar{y}_{ij}^c(s)\}},\sum_{i,j=B_n}^J m_{ij}^c(s)\delta_{\{x_{ij}^c(s),y_{ij}^c(s)\}}\right)}=\mathcal{O}(h^2)$

Using that the characteristics of couples satisfy $(x_{ij}^c(s),y_{ij}^c(s))=(\bar{x}_{ij}^c(s),\bar{y}_{ij}^c(s))$ and the evolution of the masses $m_{ij}^c(s)$ and $\bar m_{ij}^c(s)$, we get that
\begin{align*}
\sum_{i,j=B_n}^J \left|m_{ij}^c(s)-\bar m_{ij}^c(s)\right| \leq &\,\int_\tau^{\tau+h}\sum_{i,j=1}^N \left|\frac{\tilde{N}_{ij}(t)}{\tilde{D}_{ij}(t)}-\frac{N_{ij}(t)}{D_{ij}(t)}\right|dt
\\
&\hspace{-1cm}+ \int_\tau^{\tau+h}\!\!\!\sum_{i,j=1}^N \left|\xi^c(t,\nu_t^m,\nu_t^f,\nu_t^c)(x_{ij}(t),y_{ij}(t))m_{ij}(t)-
\xi^c(t,\bar{\nu}_t^m,\bar{\nu}_t^f,\bar{\nu}_t^c)(x_{ij}(t),y_{ij}(t))\bar{m}_{ij}(t)\right| dt\\
\leq & \,h^2C_9+h^2\|\xi^c\|_{\mathbf{BC}}(2C_2C_3+C_2)
\end{align*}
where $C_9=\left(8(\|h\|_{\infty,{\rm \bf Lip}}+\|g\|_{\infty,{\rm \bf Lip}})\|\Theta\|_{\infty,{\rm \bf Lip}}\right) C_2$. The estimate on $C_9$ was derived due to the observation that 
$$
\sum_{i,j=1}^N \left|\frac{\tilde{N}_{ij}(t)}{\tilde{D}_{ij}(t)}-\frac{N_{ij}(t)}{D_{ij}(t)}\right| \leq d_2\left(\mathcal{T}(\nu_t^f,\nu_t^m,\nu_t^c),\mathcal{T}(\bar{\nu}_t^f,\bar{\nu}_t^m,\bar{\nu}_t^c)\right),
$$
taking into account \eqref{spec_marriage_function_measures}. Moreover, using \cite[Lemma 2.4]{Ulikowska2012}, one obtains
\begin{align*}
d_2\left(\mathcal{T}(\nu_t^f,\nu_t^m,\nu_t^c),\mathcal{T}(\bar{\nu}_t^f,\bar{\nu}_t^m,\bar{\nu}_t^c)\right)& \leq\|\mathcal{T}\|_{{\rm \bf BC}^{0,1}}\left( d_1(\nu_t^f,\bar{\nu}_t^f)+d_1(\nu_t^m,\bar{\nu}_t^m)+d_2(\nu_t^c,\bar{\nu}_t^c)\right)
\\
&\leq C_2 h \left(8(\|h\|_{\infty,{\rm \bf Lip}}+\|g\|_{\infty,{\rm \bf Lip}})\|\Theta\|_{\infty,{\rm \bf Lip}}\right) \,.
\end{align*}
\end{enumerate}
\smallskip
Obviously distance $d_1(\bar{\nu}_{s}^f,\nu_{s}^f)$ can be estimated analogously to $d_1(\bar{\nu}_{s}^m,\nu_{s}^m)$.
It just has been shown that 
$$
{\mathbf d}({\textbf v}_{\tau+h},S(h;\tau){\textbf v}_\tau)=(\mathcal{O}(h^2)+\mathcal{O}(h))\Delta t+\mathcal{O}(h^2).
$$
Applying this local estimate to Proposition \ref{tangential_ineq} we obtain the claim of the underlying theorem, namely
$$
\mathbf{d}(\mathbf{\mathbf{u}_t},\mathbf{\mathbf{v}_t})\leq 
\varepsilon_0+
L\int_{[0,t]}\liminf_{h\rightarrow 0}\frac{{\mathbf d}({\textbf v}_{\tau+h},S(h;\tau){\textbf v}_\tau)}{h}d
\tau\leq
\varepsilon_0+L\int_{[0,t]} C \Delta t\ d\tau
\leq\varepsilon_0+LCt\Delta t.
$$
The constant $C$ estimates a certain combination of constants $C_1,\ldots,C_7$.
\end{proofof}

\begin{Remark}\label{internalisation_point}
At the beginning of the proof we assumed, that interval $(\tau,\tau+h]$ does not contain the internalisation point $t_n$, then we consider $(\tau,\tau+h]=(t_n,\tau+h]$, where $\tau+h<t_{n+1}$ and observe, that $\bar{m}_{B_n}^m=\bar{m}_{B_n}^f=0$, that $\bar{m}_{ij}^c(t)=m_{ij}^c(t)=0$, $j=B_n \lor\ j=B_n$, and that the whole argumentations of the proof does not change.
\end{Remark}


\section{Numerical examples}
\setcounter{equation}{0}

In this Section we present two numerical examples, illustrating the theoretical results. In both cases we present tables of errors indicating the rate of the convergence. In the second example we calculate the errors not only in bounded Lipschits distance, but also in TV. As it turns our the convergence cannot be obtained in TV, what confirms the need of presented theory. The measurement of the error in bounded Lipschits distance is truly necessary in those calculations, but is far from trivial and requires additional explanation. For this reason we start with the details concerning the error measurement in Subsection \ref{measure_of_error}. Subsections \ref{example_1} and \ref{example_2} deal with the numerical examples. In Subsection \ref{example_1} we consider an equation with the simplest possible coefficients satisfying the theoretical assumptions, whose solution is not known. For this reason Table \ref{tab:table1} presents the errors between the numerical solutions and a reference solution. It is also worth to notice here, that due to the choice of trivial/constant mortality rates ($c^m=c^f=0.1$), the number of males and females within any cohort never reaches $0$, what is a non-realistic phenomena, because usually it is assumed that every single individual is eventually dying. To present some more probable model, we introduce the second example in Subsection \ref{example_2}, where the the coefficients are somewhat complicated and time-dependent. They were chosen not only to satisfy the theoretical assumptions, but also in such a way,  that we know the exact solution of the system. More precisely we first imposed the equations describing the evolution of males, females and couples, and derived the coefficient which meet the assumptions and for which the imposed equations satisfy the system. Given the exact solution we can observe that all the individuals, as well as all the couples in any cohort, will eventually extinct. This is in accordance to life observations. Moreover, knowing the exact solution of the system, the errors presented in Table \ref{tab:table2} measure the distance between the numerical and the analytical solutions.

\subsection{Measurement of the error}\label{measure_of_error}

Due to the definition of bounded Lipschitz distance (which is a supremum over bounded Lipschitz functions), the calculation of an error in flat metric is not straightforward. First of all we need to consider approximation of initial data and its error. Secondly we wish to reduce the problem of computational bounded Lipschitz distance between two atomic measures to the problem of computational $1$-Wasserstein distance. Thirdly the computational cost of Wasserstein distance in higher dimension (2D) is troublesome itself and deserves a special attention.

Every numerical computation starts with establishing the time and space steps, here $\Delta t$, $\Delta x$ and $\Delta y$. Given that in the case of aged structured population models characteristic are straight lines, it is natural to assume that $\Delta t=\Delta x=\Delta y$ and to divide the domain of solution into $T/\Delta t$ time steps,  $M/\Delta t$ space cells in case of male and female population and $(M/\Delta t)^2$ in case of couple population with $M$ being the largest age of the population. According to Remark  \ref{initial_accuracy}, the initial condition can approximated arbitrarily accurate by a linear combination of certain amount of Dirac measures. Special techniques of measure reconstruction described in \cite{MR3244779} allow to present this desired approximation with a linear combinations of only $M/ \Delta x$ Dirac Deltas (or $(M/ \Delta x)^2$ in case of couples) with error $\Delta x^2$. It is easy to notice that the same reasoning can be applied in 2D case. Given an exact solution at time $T$ in Subsection \ref{example_2}, we can use the same methods. Analysing the procedures of measure reconstruction proposed in \cite{MR3244779},  it is sufficient to approximate initial conditions with formulas (\ref{the_numerical_method_initial_1a}) and (\ref{the_numerical_method_initial_1b}), to attain satisfactory approximation $\varepsilon_0={\mathcal O}(\Delta t^2)$.
Obviously, presenting a method of order one only, this inaccuracy can be neglected.

Lemma 2.1 in \cite{MR3244779} shows how to reduce the problem of bounded Lipschitz distance (in 1D) to some other measure expressed in terms of $1$-Wasserstein distance.The same reasoning can be easily adapted to two dimensional case. 

\begin{Lemma}
Let $\mu_1,\mu_2\in\mathcal{M}_+(\bR_+)$ be such that $M_{\mu_i}=\int_{\bR_+}d\mu_i\neq0$ and $\tilde{\mu}_i=\mu_i/M_{\mu_i}$ for $i=1,2$. Define $\rho:\mathcal{M}_+(\bR_+)\times \mathcal{M}_+(\bR_+)\rightarrow\bR_+$ as the following
$$
\rho(\mu_1,\mu_2)={\rm min}\{M_{\mu_1},M_{\mu_2}\}W_1(\tilde{\mu}_1,\tilde{\mu}_2)+|M_{\mu_1}-M_{\mu_2}|,
$$
where $W_1$ is the $1$-Wasserstein distance. Then, there exists a constant $C_K=\frac{1}{3}{\rm min}\left\{1,\frac{2}{|K|}\right\}$, such that
$$
C_K\rho(\mu_1,\mu_2)\leq d_1(\mu_1,\mu_2)\leq\rho(\mu_1,\mu_2),
$$
where K is the smallest interval such that ${\rm supp}(\mu_1)$, ${\rm supp}(\mu_2)\subseteq K$ is the length of the interval $K$. If $K$ is unbounded we set $C_K=0$.
\end{Lemma}

In all presented numerical experiments, the effective error of the method ${\rm Err}(\Delta t)$ will be estimated in terms of metric $\rho$.
To compute effectively the Wasserstein distance $W_1(\tilde{\mu}_1,\tilde{\mu}_2)$ in any dimension, we resort to the results presented in \cite{BCCNP15}, where the considerations start from approximation of  $\tilde{\mu}_1$, $\tilde{\mu}_2$ by some atomic measures $\sum_i^{N_a} a_i \delta_{\{ x^a_i\}}$ and $\sum_j^{N_b} b_j \delta_{\{ x^b_j\}}$ respectively (for sake of simplicity we assume that $N_a=N_b$). Instead of computing
\begin{equation}\label{wass_atomic}
W_1\left(\sum_i a_i \delta_{\{ x^a_i\}},\sum_j b_j \delta_{\{ x^b_j\}}\right):=\min\left\{ \sum_{i,j} \left(c_{ij}\gamma_{ij}\right):\quad \gamma_{ij}\geq0,\ \sum_i\gamma_{ij}=b_j,\  \sum_j\gamma_{ij}=a_i\right\}
\end{equation}
we fix small $\varepsilon>0$ and focus on
\begin{equation*}
W_1^{\varepsilon}:=\min\left\{ \sum_{i,j} \left(c_{ij}\gamma_{ij}+\varepsilon\gamma_{ij}\log(\gamma_{ij})\right):\quad \gamma_{ij}\geq0,\ \sum_i\gamma_{ij}=b_j,\  \sum_j\gamma_{ij}=a_i\right\},
\end{equation*}
which for $\varepsilon\rightarrow0$ tends to the minimization problem (\ref{wass_atomic}) in the sense of $\Gamma$-convergence. Taking $\eta_{ij}:={\rm e}^{-c_{ij}/\varepsilon}$, we observe that
$$
c_{ij}\gamma_{ij}+\varepsilon\gamma_{ij}\log(\gamma_{ij})=\varepsilon\gamma_{ij}\log\left(\frac{\gamma_{ij}}{\eta_{ij}}\right)=\varepsilon{\rm KL}(\gamma|\eta),
$$
where ${\rm KL}$ is  Kullback-Leiber divergence, that is a sort of a distance based on a relative entropy:
$$
{\rm KL}(\gamma|\eta):=
\left\{
\begin{array}{ll}
\sum_{ij}\gamma_{ij}\log\left(\frac{\gamma_{ij}}{\eta_{ij}}\right)\ & {\rm if}\ \frac{\gamma_{ij}}{\eta_{ij}}>0,\\
0\ & {\rm if}\ \frac{\gamma_{ij}}{\eta_{ij}}=0,\\
+\infty,\ & {\rm if}\ \frac{\gamma_{ij}}{\eta_{ij}}<0,\\
\end{array}
\right.\\
$$
Given a convex set $\mathcal{C}\in R^{N_a\times N_a}$, the projection according to the Kullback-Leiber divergence is defined as
$$
P_{\mathcal{C}}^{\rm KL}(\eta):={\rm argmin}_{\gamma\in\mathcal{C}}{\rm KL}(\gamma|\eta).
$$
This means that
$
W_1^{\varepsilon}={\rm KL}(P_{\mathcal{C}}^{\rm KL}(\eta)|\eta),
$
where $P_{\mathcal{C}}^{\rm KL}(\eta)$ can be computed using {\em Iterative Bergman Projections}:
$$
\gamma^{(0)}:={\rm e}^{-C/\varepsilon},\ \gamma^{(n)}:=P_{\mathcal{C}}^{\rm KL}( \gamma^{(n-1)}),
$$
with the entries of $C$ defined as $c_{ij}=\|x_i^a-x_j^b\|$.
It can be shown that
$$
\gamma^{(n)} \rightarrow P_{\mathcal{C}}^{\rm KL}(\eta),\  {\rm as}\  n\rightarrow \infty.
$$
For more details we refer to \cite{OTAM}. The rate of convergence $q$ presented in the tables of errors is given by
$$
q:=\lim_{\Delta t\rightarrow0}\frac{\log[{\rm Err}(2\Delta t)/{\rm Err}(\Delta t)]}{\log2}.
$$

\subsection{Example 1}\label{example_1}
In the first example we approximate system (\ref{the_model_densities_nonlinear}) for  $t\in[0,1)$ and $(x,y)\in[0,1)\times[0,1)$, where mortality and birth rates are constant
$$
c^m(t,x)=c^f(t,y)=c^c(t,x,y)=0.1
$$
$$
b^m(t,x,y)=b^f(t,x,y)=10,
$$
and marriage function coefficients $h$, $g$ and $\Theta$ do not depend on time
$$
h(t,x)=\left\{
\begin{array}{cc}
 \left(\frac{1}{10}-x\right) (x-1) & \frac{1}{10}\leq x\leq 1 \\[2mm]
 0 & \text{otherwise} \\
\end{array}\right.
\quad \quad
g(t,y)=\left\{
\begin{array}{cc}
 \left(\frac{1}{10}-y\right) (y-1) & \frac{1}{10}\leq y\leq 1 \\[2mm]
 0 & \text{otherwise} \\
\end{array}\right.
$$
$$
\Theta(t,x,y)=\left\{
\begin{array}{cc}
10 \left(\frac{1}{10}-x\right) (x-1) \left(\frac{1}{10}-y\right) (y-1) & \frac{1}{10}\leq x\leq 1,\frac{1}{10}\leq y\leq 1 \\[2mm]
 0 & \text{otherwise} \\
\end{array}\right.\,.
$$
The initial conditions for the system are also trivial and given by
$$
u^m(0,x)=u^f(0,y)=u^c(0,x,y)=1.
$$
Table \ref{tab:table1} shows that the scheme behaves as an scheme of order 1 as expected from our Theorem \ref{convergence}.

\begin{table}[h!]
  \centering
  \begin{tabular}{rrr}
    \toprule
    $\Delta t=\Delta x$ & ${\rm Err}(\Delta t)$ & $q$ \\
    \midrule
    $10^{-1}$ &  $5.89\cdot10^{-2}$ & $-$\\
    $5*10^{-2}$ & $2.57\cdot 10^{-2}$ & $1.196499275$\\
    $2.5*10^{-2}$ & $1.2\cdot 10^{-2}$ & $1.09873395$\\
    $1.25*10^{-2}$ & $5.9\cdot 10^{-3}$ & $1.02424755$\\
    $6.25*10^{-3}$ & $2.93 \cdot 10^{-3}$ & $1.00981429$\\
    $3.125*10^{-3}$ & $1.46 \cdot 10^{-3}$ & $1.0049323$\\
    $1.5625*10^{-3}$ & $7.39 \cdot 10^{-4}$ & $0.9823221$\\
     $7.8125*10^{-4}$ &$3.7 \cdot 10^{-4}$ & $0.99804909$\\
    \bottomrule
  \end{tabular}
   \caption{Error computed in flat metric ${\rm Err}(\Delta t)$ and its order of convergence $q$.} \label{tab:table1}
\end{table}

\subsection{Example 2}\label{example_2}

We now present a numerical example for the system (\ref{the_model_densities_nonlinear}), whose exact solution is known, evolves in $[0,\infty)^3$ and is given by the following formulas
$$
u^m(t,x)=\left\{
\begin{array}{cc}
 \left(1-\frac{t}{10}\right) (t-x-1) (-t+x-1) & 0\leq x\leq t+1 \\[2mm]
 0 & \text{otherwise} \\
\end{array}
\right.
$$
$$
u^f(t,y)=\left\{
\begin{array}{cc}
 \left(1-\frac{t}{10}\right) (t-y-1) (-t+y-1) & 0\leq y\leq t+1 \\[2mm]
 0 & \text{otherwise} \\
\end{array}
\right.
$$
$$
u^c(t,x,y)=\left\{
\begin{array}{cc}
\left(1-\frac{t}{10}\right) \left(x-\frac{1}{10}\right)^2 \left(y-\frac{1}{10}\right)^2 (-t+x-1)^2 (-t+y-1)^2 & \frac{1}{10}\leq x\leq t+1\land \frac{1}{10}\leq y\leq t+1 \\[2mm]
 0 & \text{otherwise} \\
\end{array}
\right.
$$
The advantage of having the exact solution goes in pair with the disadvantage of lengthy and complicated coefficients given by
$$
c^m(t,x)=\left\{
\begin{array}{cc}
 \frac{1}{10-t} & 0\leq x\leq t+1 \\[2mm]
 0 & \text{otherwise} \\
\end{array}
\right.\,,
\quad \quad
c^f(t,y)=\left\{
\begin{array}{cc}
 \frac{1}{10-t} & 0\leq y\leq t+1 \\[2mm]
 0 & \text{otherwise} \\
\end{array}
\right.\,,
\quad \quad
c^c(t,x,y)=0.1\,,
$$
$$
b^m(t,x,y)=b^f(t,x,y)=\left\{
\begin{array}{cc}
-\frac{9000000000000 (t-1) (t+1)}{(10 t+9)^{10}} & 10 x\geq 1\land t+1\geq x\land 10 y\geq 1\land t+1\geq y \\[2mm]
 0 & \text{otherwise} \\
\end{array}
\right.\,,
$$
$$
h(t,x)=\left\{
\begin{array}{cc}
 \left(\frac{1}{10}-x\right) (-t+x-1) & \frac{1}{10}\leq x\leq t+1 \\[2mm]
 0 & \text{otherwise} \\
\end{array}\right.\,,
\quad \quad
g(t,y)=\left\{
\begin{array}{cc}
 \left(\frac{1}{10}-y\right) (-t+y-1) & \frac{1}{10}\leq y\leq t+1 \\
 0 & \text{otherwise} \\
\end{array}\right.\,,
$$
and
$$
\Theta(t,x,y)=\frac{{\rm num}(t,x,y)}{\rm den(t)}
$$
where
{\tiny
$$
{\rm num}(t,x,y)=
\frac{-\frac{t (10 x (10 y+199)+1990 y-399)}{10000}+2 x+2 y-\frac{2}{5}}
{\left(\frac{(t-10) (10 t+9)^5 (1-10 x)^2 (-t+x-1)}{3000000000}+\left(1-\frac{t}{10}\right) (t-x-1)\right) \left(\frac{(t-10) (10 t+9)^5 (1-10 y)^2 (-t+y-1)}{3000000000}+\left(1-\frac{t}{10}\right) (t-y-1) \right)}
$$}
and {\tiny
$$
{\rm den}(t)=
\frac{(t-10) (10 t+9)^4 \left(10^8 t^8+72\cdot10^7 t^7+2268\cdot 10^6 t^6+40824\cdot 10^5 t^5+45927\cdot 10^5 t^4+3306744\cdot 10^3 t^3+1488034800 t^2+21382637520 t-51056953279\right)}{21\cdot 10^{15}}+1\,.
$$}
Given the above coefficients $h$, $g$ and $\Theta$ one can check that the marriage function is given with
\begin{align*}
T(t,x,y)=&\,\frac{1}{(-t+x-1) (-t+y-1)}
\frac{(t-10) (10 t+9)^5 (1-10 x)^2 (t-x+1)^2}{3000000000}\\
&+\left(1-\frac{t}{10}\right) (t-x-1) (-t+x-1)
\frac{(t-10) (10 t+9)^5 (1-10 y)^2 (t-y+1)^2}{3000000000}\\
&+\left(1-\frac{t}{10}\right) (t-y-1) (-t+y-1)\,.
\end{align*}
We first measured the error in TV, and later in flat metric. According to our expectations, TV does not show any convergence, see Table \ref{tab:table2}, while flat metric significantly decreases the error. The rate of the error in Table \ref{tab:table2} is 1 as expect from the theoretical result in Theorem \ref{convergence}.

\begin{table}[h!]
  \centering
  \begin{tabular}{crrr}
    \toprule
    ${\rm Err_{TV}}(\Delta t)$  & $\Delta t=\Delta x$ & ${\rm Err}(\Delta t)$  & $q$ \\
    \midrule
    $8.12*10^{-2}$ & $10^{-1}$ &  $7.16\cdot 10^{-2}$ & $-$ \\
    $7.04*10^{-2}$ & $5*10^{-2}$ & $4.0\cdot 10^{-2}$ & $0.839959587$\\
    $6.89*10^{-2}$ & $2.5*10^{-2}$ & $2.13\cdot 10^{-2}$ & $0,90914657$\\
    $6.72*10^{-2}$ & $1.25*10^{-2}$ & $1.14\cdot 10^{-2}$ & $0,901819606$\\
    $6.72*10^{-2}$ & $6.25*10^{-3}$ & $6.3\cdot 10^{-3}$ & $0,855610091$\\
    $6.71*10^{-2}$ & $3.125*10^{-3}$ & $3.2\cdot 10^{-3}$ & $0,977279923$\\
    $6.70*10^{-2}$ & $1.5625*10^{-3}$ & $1.6\cdot 10^{-3}$ & $1$\\
    $6.70*10^{-2}$ & $7.8125*10^{-4}$ &$8.0\cdot 10^{-4}$ & $1$\\
    \bottomrule
  \end{tabular}
   \caption{Error computed in TV, ${\rm Err_{TV}}(\Delta t)$,  error computed in flat metric, ${\rm Err}(\Delta t)$, and order of convergence in obtained in flat metric $q$.}
   \label{tab:table2}
\end{table}

\section*{Acknowledgments}
JAC was partially supported by the EPSRC grant number EP/P031587/1.
PG was financed by The National Center for Science UMO-2015/18/M/ST1/00075.
KK was financed by The National Center for Science DEC-2012/05/E/ST1/02218.
AMC was supported by the Emmy Noether Programme of the German Research Council (DFG).
This work was partially supported by the grant 346300 for IMPAN from the Simons Foundation and the matching 2015-2019 Polish MNiSW fund.\\

\bibliographystyle{abbrv}
\bibliography{EBT-Convergence-Age_struct_two_sex_pop_model}

\end{document}